# STABILIZATION OF THE WAVE EQUATION ON LARGER-DIMENSION TORI WITH ROUGH DAMPINGS

## M. ROUVEYROL

ABSTRACT. This paper deals with uniform stabilization of the damped wave equation. When the manifold is compact and the damping is continuous, the geometric control condition is known to be necessary and sufficient. In the case where the damping is a sum of characteristic functions of polygons on a two-dimensional torus, a result by Burq-Gérard states that stabilization occurs if and only if every geodesic intersects the interior of the damped region or razes damped polygons on both sides. We give a natural generalization of their result to a sufficient condition on tori of any dimension $d \geq 3$. In some particular cases, we show that this sufficient condition can be weakened.

## Contents

1. NOTATION AND MAIN RESULT

Consider the damped wave (or Klein-Gordon) equation

(1.1) $\quad (\partial_t^2 - \Delta_g + a(x)\partial_t + m(x))u = 0, (u|_{t=0}, \partial_t u|_{t=0}) = (u_0, u_1) \in H^1(M) \times L^2(M),$

where $(M, g)$ denotes a smooth compact Riemannian manifold without boundary, $g$ is the manifold metric, $\Delta_g$ is the Laplace operator on $M$, and the damping $a$ and potential $m$ are two non-negative $L^\infty$ functions over $M$. The energy

(1.2) $$E_m(u)(t) = \frac{1}{2} \int_M (|\nabla_g u|_g^2 + |\partial_t u|^2 + m|u|^2) d\operatorname{vol}_g$$

is then decaying, as

$$E_m(u)(t) = E_m(u)(0) - \int_0^t \int_M a(x)|\partial_t u(t,x)|^2 d\operatorname{vol}_g(x) dt.$$

We say that *uniform stabilization* holds for the damping $a$ if one of the following equivalent conditions is satisfied:

(1) There exists a rate $f(t)$ such that $\lim_{t \to +\infty} f(t) = 0$ and for any $(u_0, u_1) \in H^1(M) \times L^2(M)$,
$$E_m(u)(t) \le f(t) E_m(u)(0).$$

(2) There exist some constants $C, c > 0$ such that for any $(u_0, u_1) \in H^1(M) \times L^2(M)$,
$$E_m(u)(t) \le Ce^{-ct} E_m(u)(0).$$

(3) There exist some $T > 0$ and $c > 0$ such that for any $(u_0, u_1) \in H^1(M) \times L^2(M)$, if $u$ is a solution of the damped wave equation (1.1), then
$$E_m(u)(0) \le C \int_0^T \int_M 2a(x)|\partial_t u|^2 d\operatorname{vol}_g.$$

(4) There exist some $T > 0$ and $c > 0$ such that for any $(u_0, u_1) \in H^1(M) \times L^2(M)$, if $u$ is a solution of the *un*damped wave equation
$$(\partial_t^2 - \Delta + m)u = 0, \quad (u|_{t=0}, \partial_t u|_{t=0}) = (u_0, u_1) \in H^1(M) \times L^2(M),$$
then
$$E_m(u)(0) \le C \int_0^T \int_M 2a(x)|\partial_t u|^2 d\operatorname{vol}_g.$$

The HUM method developed by J.-L. Lions [Lio88] (see also [BG02, Chapter 4.2]) provides a functional-analytic framework which reduces controllability and stabilization problems to observability inequalities and unique continuation properties for a dual equation. In the case of the wave equation, it allows to prove that uniform stabilization is equivalent to observability estimates (see Proposition 2.1 for a precise statement) and to exact controllability of the equation. We refer to [BG02, Chapters 3 and 4] for proofs concerning the wave equation and its controllability properties.

In the case where the damping $a$ is continuous, a necessary and sufficient condition for uniform stabilization is given by the following landmark result:

**Theorem 1** (Geometric control condition, [RT75, BLR92, BG97]). — *Let $m \ge 0$. Assume that the damping $a$ is continuous. For $\rho_0 = (x_0, \xi_0) \in S^*M$ (the unit cotangent bundle of* $M$), let $\gamma_{\rho_0}(s)$ denote the geodesic starting from $x_0$ in the (co-)direction $\xi_0$. The damping $a$ then stabilizes uniformly the wave equation if and only if:

(GCC) $\quad \exists T, c > 0$ such that $\inf_{\rho_0 \in S^*M} \int_0^T a(\gamma_{\rho_0}(s))ds \geq c,$

or equivalently if every geodesic $\gamma_{\rho_0}(s)$ intersects $\{a > 0\}$ in time $T$.

In the case where $a$ is merely $L^\infty$, the following classical result is a consequence of Bardos, Lebeau and Rauch's work:

**Theorem 2** ([BLR92], [BG20]-Theorem 2). — *Assume that $0 \leq a \in L^\infty(M)$. Then the strong geometric control condition*

(SGCC) $\quad \begin{aligned} &\exists T, c > 0 \text{ such that } \forall \rho_0 \in S^*M, \exists s \in (0,T), \delta > 0 \\ &\text{such that } a \geq c \text{ a.e. over } B(\gamma_{\rho_0}(s), \delta) \end{aligned}$

*is **sufficient** for uniform stabilization, and the weak geometric control condition*

(WGCC) $\quad \exists T > 0$ *such that* $\forall \rho_0 \in S^*M, \exists s \in (0,T)$ *such that* $\gamma_{\rho_0}(s) \in \text{supp}(a),$

*where $\text{supp}(a)$ is the support of $a$ in the distributional sense, is **necessary** for uniform stabilization.*

Finding a necessary and sufficient condition in between (SGCC) and (WGCC) in the case of $L^\infty$ dampings yet remains an open problem. In [Leb92], Gilles Lebeau proved that uniform stabilization holds on the sphere $\mathbb{S}^d = \{x \in \mathbb{R}^{d+1}, |x| = 1\}$ when $a = 1$ over the half-sphere. Hui Zhu generalized this result to Zoll surfaces of revolution in [Zhu18]. See also [AL14, Appendix B] for a case where the damped region $\{a = 1\}$ does not satisfy (WGCC).

A necessary and sufficient condition for uniform stabilization was given by Nicolas Burq and Patrick Gérard in [BG20], in the case where $M$ is a 2-dimensional torus $\mathbb{T}^2 = \mathbb{R}^2 / A\mathbb{Z} \times B\mathbb{Z}$ with $A, B > 0$, and the damping $a$ is a sum of characteristic functions of disjoint polygons. Their result is stated thereafter and illustrated in Figure 1.1.

**Theorem 3** ([BG20], Assumption 1.2 and Theorem 4). — *The damping $a$ stabilizes the wave equation if and only if the following assumption holds: there exists some $T > 0$ such that all geodesics (straight lines) of length $T$ either encounter the interior of one of the polygons or follow for some time one of the sides of a polygon on the left and for some time one of the sides of a polygon (possibly the same) on the right.*

The present article is dedicated to exploring generalizations of Theorem 3 to higher-dimensional tori. Like in [BG20], we use dampings $a$ that are characteristic functions. This is motivated both by the fact that applied control historically uses characteristic functions, and by a will to understand the concentration and non-concentration properties of waves in situations where geodesics raze the control region without entering its interior. The case of tori is a favorable one for this study because the absence of curvature or boundary simplifies the geometry, and because it ensures that quasimodes of the Laplace operator satisfy good non-concentration properties (see Proposition 3.1).

We identify a $d$-dimensional torus $\mathbb{T}^d$ with a cuboid $\prod_{j=1}^d [0, A_j]$, $A_j > 0$, with appropriate identification of the faces. We call a *polyhedron* of $\mathbb{T}^d$ the intersection of the cuboid with a finite number of half-spaces, and we assume that $a$ is a finite sum of characteristic

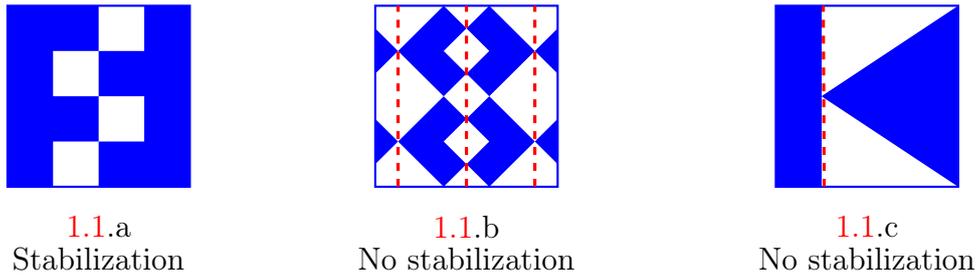

1.1.a
Stabilization

1.1.b
No stabilization

1.1.c
No stabilization

FIGURE 1.1. N. Burq and P. Gérard's checkerboards on $\mathbb{T}^2$ [BG20]: the damping $a$ is equal to 1 in the colored region and 0 elsewhere. For all these examples, (WGCC) is satisfied but not (SGCC). The dashed lines are geodesics which violate the condition of Theorem 3.

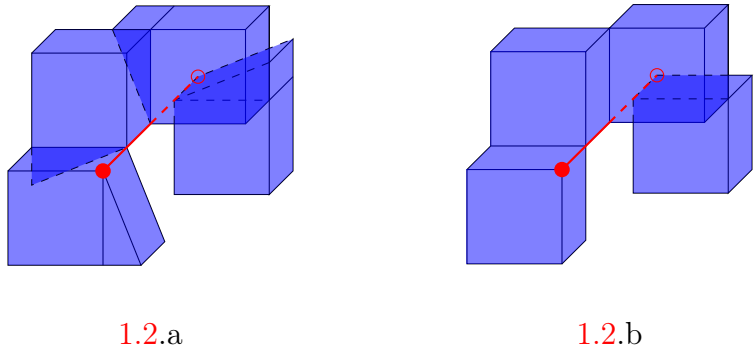

1.2.a
1.2.b

FIGURE 1.2. In these three-dimensional illustrations, the damping $a$ is equal to 1 in the colored zone and 0 elsewhere. Figure 1.2.a illustrates the framework of Theorem 4, as the red geodesic satisfies condition (1.3) but does not intersect the interior of supp($a$). In Figure 1.2.b, the red geodesic does not satisfy this condition as directions $\Xi = (\pm 1, 0)$ or $(0, \pm 1)$ are not damped in the sense of Definition 1.1. In Section 5, we prove that uniform stabilization can still occur in cases involving such geodesics.

functions of disjoint such polyhedrons. This choice comes from the fact that it is always possible to concentrate solutions of the wave equation near a geodesic that only has finite-order contacts with the damped region (see for example [BG20, Section 5]). Thus, geodesics which violate (SGCC) need to have infinite-order contacts with the damped zone for stabilization to occur. By considering a damping equal to 1 over a finite union of polyhedrons, we essentially cover all cases where the damping is a characteristic function.

On $\mathbb{T}^2$, the conormal space to a geodesic at a point is one-dimensional, hence the unit sphere of (co-)normal directions to the geodesic is reduced to two points which correspond to the right and left sides of the geodesic. The assumption of Theorem 3 thus naturally generalizes to assuming that every geodesic is damped in each of its normal directions. Our first result states that this condition is indeed sufficient for uniform stabilization. It holds on any torus $\mathbb{T}^d = \mathbb{R}^d/(A_1\mathbb{Z} \times \cdots \times A_d\mathbb{Z})$ of dimension $d \geq 3$.

**Theorem 4.** — *Assume that every geodesic orbit of $\mathbb{T}^d$ either intersects the interior of* $\mathrm{supp}(a)$, *or is damped in every normal direction, ie*

(1.3)
$$\text{for all } \rho_0 = (X_0, \Xi_0) \in S^*\mathbb{T}^d, \Xi \text{ orthogonal to } \Xi_0 \text{ with } |\Xi| = 1,$$
$$\text{there exist some positive } \delta_0 \text{ and some interval } I \subset \mathbb{R} \text{ such that}$$
$$\gamma_{\rho_0}(s) + \delta\Xi \in \mathrm{Int}(\mathrm{supp}(a)), \quad \forall s \in I, \quad \forall \delta \in (0, \delta_0].$$

*Then stabilization holds.*

Note that any geodesic which encounters the interior of $\mathrm{supp}(a)$ also satisfies assumption (1.3). A typical example of a geodesic satisfying the assumption without entering the interior of $\mathrm{supp}(a)$ is the red one in Figure 1.2.a. The proof of Theorem 4 consists in studying precisely the concentration of high-frequency quasimodes of the Laplacean near such a geodesic (specifically near the red geodesic of Figure 4.1, see section 4.2) using second microlocalization techniques, then reducing the general case to this model case.

Condition (1.3) naturally leads to the following definition:

**Definition 1.1.** — *Consider a geodesic $\gamma_{\rho_0}$ of $\mathbb{T}^d$ with direction $\Xi_0 \in \mathbb{S}^{d-1}$, and some $\Xi \in \mathbb{S}^{d-1}$ orthogonal to $\Xi_0$. We say that $\gamma_{\rho_0}$ is damped in the normal direction $\Xi$ if there exists some positive $\delta_0$ and some interval $I \subset \mathbb{R}$ such that*

(1.4) $$\gamma_{\rho_0}(s) + \delta\Xi \in \mathrm{Int}(\mathrm{supp}(a)), \quad \forall s \in I, \quad \forall \delta \in (0, \delta_0].$$

For example, in Figure 1.2.b, the red geodesic is damped in all but four normal directions. Note also that since the interval $I$ is independent of $\delta$ in the definition, a geodesic can be damped in some direction only if it razes an edge or a face of some polyhedron in that direction. Punctual contacts between a damped polyhedron and a geodesic are not considered sufficient for the geodesic to be damped in the normal directions that enter the polyhedron, so as to avoid phenomena resembling that of Figure 1.1.c.

We emphasize that in dimension 2, assumption (1.3) is exactly the necessary and sufficient condition found by N. Burq and P. Gérard in [BG20]. In dimensions 3 and greater, the enriched geometry (specifically the fact that the unit sphere of the conormal space to a geodesic at a point is infinite) permits uniform stabilization under weaker conditions, as in the following result:

**Theorem 5.** — *Identify $\mathbb{T}^3 \simeq [-1,1]_x^2 \times [-1,1]_y$ and assume that $a$ is the characteristic function of the set*

(1.5)
$$\{|x_1| \geq \tfrac{1}{2} \text{ or } |x_2| \geq \tfrac{1}{2}\} \cup \{-1 < y < -\tfrac{1}{2}, x_1 > 0, x_2 < 0\} \cup \{-\tfrac{1}{2} < y < 0, x_1 > 0, x_2 > 0\}$$
$$\cup \{0 < y < \tfrac{1}{2}, x_1 < 0, x_2 > 0\} \cup \{\tfrac{1}{2} < y < 1, x_1 < 0, x_2 < 0\},$$

*as illustrated in Figure 5.1. Then stabilization holds.*

This case is based on the shape of the damping in Figure 1.2.b, with some geodesics damped in all but a finite number of directions, so that condition (1.3) is not satisfied. We prove this theorem using directional second microlocalization and third microlocalization procedures that are new to the best of our knowledge.

The paper is constructed as follows. In Section 2, we recall the first microlocalization procedure for the wave equation and introduce relevant tools. In Section 3, we prove non-concentration estimates for Laplacean quasimodes on $\mathbb{T}^d$. Section 4 is dedicated to the proof of Theorem 4 and Section 5 to the refined example of Theorem 5. In the concluding paragraphs of Section 5, we discuss possible extensions of this theorem using our techniques.

*Acknowledgments* The author wishes to thank Nicolas Burq for his rich and careful advice during the writing of this article and Antoine Prouff for the many discussions about it. We also thank an anonymous referee for the time and attention they dedicated to reviewing this paper. This work was funded by a CDSN PhD grant from Ecole Normale Supérieure Paris-Saclay via the Hadamard Doctoral School of Mathematics (EDMH).

## 2. First microlocalization

The goal of this section is to introduce the first microlocal measure and the observability estimate that will be used later on. We do so by proving that (SGCC) implies uniform stabilization. The result of this section is valid on any compact Riemannian manifold without boundary $M$ and for any non-negative damping function $a \in L^\infty(M)$. For simplicity, we prove it on $d$-dimensional tori. We refer to [BG20, Sections 2, 5 and Annex A] for statements and proof corresponding to (SGCC) and (WGCC) in the general case.

Let us recall the statement of the strong geometric control condition:

**Theorem 6.** — *Assume that* $0 \leq a \in L^\infty(\mathbb{T}^d)$. *Then the* **strong geometric control condition**

$$\text{(SGCC)} \qquad \begin{aligned} &\exists T, c > 0 \text{ such that } \forall \rho_0 \in S^*\mathbb{T}^d, \exists s \in (0, T), \delta > 0 \\ &\text{such that } a \geq c \text{ a.e. over } B(\gamma_{\rho_0}(s), \delta) \end{aligned}$$

*is* **sufficient** *for uniform stabilization.*

*Proof.* Let us assume that (SGCC) holds. The following result states that uniform stabilization is equivalent to an observability estimate for solutions of the Helmholtz equation in the high-frequency regime.

**Proposition 2.1** ([BG20], Proposition A.5). — *Consider an $L^\infty$ non-negative function $a$ such that $\int_{\mathbb{T}^d} a(x)dx > 0$, then $a$ uniformly stabilizes the wave equation (1.1) if and only if*

$$\text{(Obs)} \qquad \begin{aligned} &\exists C, h_0 > 0 \text{ such that } \forall 0 < h < h_0, \quad \forall (u, f) \in H^2(\mathbb{T}^d) \times L^2(\mathbb{T}^d), \\ &(h^2\Delta + 1)u = f, \text{ there holds } \|u\|_{L^2} \leq C\left(\|a^{\frac{1}{2}}u\|_{L^2} + \frac{1}{h}\|f\|_{L^2}\right). \end{aligned}$$

We then study microlocal measures associated to a sequence of quasimodes that violate the observability estimate to prove it by contradiction, following the idea of [Leb96, Bur02].

We assume that (Obs) does not hold and want to deduce a contradiction with (SGCC) from there. We obtain sequences $(h_n) \to 0$ and $(u_n, f_n)$ that satisfy

(2.1) $\quad (h_n^2 \Delta + 1)u_n = f_n, \quad \|u_n\|_{L^2} = 1, \quad \|a^{\frac{1}{2}}u_n\|_{L^2} = o(1)_{n \to +\infty}, \quad \|f_n\|_{L^2} = o(h_n)_{n \to +\infty}.$

Given a symbol $q \in C_c^\infty(T^*\mathbb{T}^d)$, we define its quantization as follows. Using a partition of unity, we can assume that $q$ is supported in a local chart. In this chart, we define

$$\text{(2.2)} \qquad \text{Op}_h(q)u(X) = \frac{1}{(2\pi)^d} \int_{\mathbb{R}^d \times \mathbb{R}^d} e^{i(X-Y)\cdot\Xi} q(X, h\Xi)\zeta(Y)u(Y)dYd\Xi$$

where $\zeta = 1$ in a neighborhood of $p_X(\text{supp}(q))$ with $p_X : T^*\mathbb{R}^d \to \mathbb{R}^d_X$. Then, up to extracting a subsequence, $u_n$ admits a semiclassical measure $\nu$ over $T^*\mathbb{T}^d$ which satisfies

$$\lim_{n \to +\infty} (\text{Op}_{h_n}(q)u_n, u_n)_{L^2} = \langle \nu, q \rangle.$$

$u_n$ is $h_n$-oscillating, since the Pancherel identity and (2.1) give

$$\int_{|\Xi| \geq \frac{R}{h_n}} |\hat{u}_n(\Xi)|^2 d\Xi \leq R^{-4} \|h_n^2 |\Xi|^2 |\hat{u}_n|\|_{L^2}^2 \leq \frac{C}{R^4} \to_{R \to +\infty} 0.$$

Thus, $\nu$ has total mass 1. Since $\text{Op}_{h_n}(-|\Xi|^2 + 1) = h_n^2 \Delta + 1$, the asymptotic

$$\|(h_n^2 \Delta + 1)u_n\|_{L^2} = o(h_n)$$

then gives that $\nu$ is supported in the characteristic set

$$\{(X, \Xi) \in T^*\mathbb{T}^d, |\Xi|^2 = 1\}$$

and that it is invariant by the bicharectistic flow:

$$\Xi.\nabla_X \nu = 0.$$

In particular, $\nu$ can be considered as a measure over the unit cotangent bundle of the torus $S^*\mathbb{T}^d = T^*\mathbb{T}^d/\mathbb{R}_+^*$. We refer to [GL93, Bur97] and [Zwo12, Chapter 5] for the construction and properties of semiclassical measures.

Denote then

$$S = \{X \in \mathbb{T}^d, \text{ there exist } \delta > 0, c > 0 \text{ such that } a \geq c \text{ a.e. on } B(X, \delta)\}.$$

Taking some $X_0 \in S$ and corresponding $\delta, c$, we get that for every symbol $q(X, \Xi)$ supported in $B(X_0, \delta)$ in the $X$ variable,

(2.3)
$$\begin{aligned}|(\text{Op}_{h_n}(q)u_n, u_n)_{L^2}| &= |(\text{Op}_{h_n}(q)u_n, \mathbb{1}_{B(X,\delta)} u_n)_{L^2}| \\ &\leq \|\text{Op}_{h_n}(q)u_n\|_{L^2} \|\mathbb{1}_{B(X,\delta)} u_n\|_{L^2} \\ &\leq \frac{1}{c} \|\text{Op}_{h_n}(q)u_n\|_{L^2} \|au_n\|_{L^2}.\end{aligned}$$

By the Calderón-Vaillancourt theorem (see for example [Zwo12, Theorem 5.1], which is also valid for the standard quantization), $\text{Op}_{h_n}(q)$ is bounded uniformly with respect to $h_n$ over $L^2$, so that $(\text{Op}_{h_n}(q)u_n)$ is a bounded sequence in $L^2$. Since $\|au_n\| = o(1)_{L^2}$, we get that $\langle \nu, q \rangle = 0$. Hence, $\nu$ vanishes in a neighborhood of every point $\rho \in S^*_X \mathbb{T}^d$ for all $X \in S$. By (SGCC), every bicharacteristic contains one such $\rho$, so that $\nu$ is identically 0. This contradicts the fact that $\nu$ has total mass 1. □

Theorem 6 can be rephrased as saying that uniform stabilization occurs when the set $\Omega = \bigcup_\omega \omega$ captures every geodesic in finite time, where the union is over all open sets $\omega$ where $a \geq c_\omega \mathbb{1}_\omega$ a.e. for some $c_\omega > 0$. In the case where $a$ is the characteristic function of an open set, $\Omega$ is exactly the interior of $\text{supp}(a)$. Thus the result will later be used in the following form, which is a direct consequence of the proof:

**Corollary 2.2.** – *Assume $a$ is the characteristic function of an open set. If a bicharacteristic contains some $\rho = (X, \Xi) \in S^*\mathbb{T}^d$ with $X \in \text{Int}(\text{supp}(a))$, then the first microlocal measure $\nu$ has no support over this bicharacteristic.*

## 3. Non-concentration estimates

In this section, we prove the non-concentration estimates that will determine the scaling of the second and third microlocalization procedures. When performing the second microlocalization procedure, these estimates allow to avoid dealing with the trace-operator-valued part of the 2-microlocal measure supported at a finite distance from the origin (see [AM14, Section 3], [FK00, Theorem 1]). On manifolds where such estimates should fail (or rather, hold for $h_n^{\frac{1}{2}}$-thickened geodesics without the $\epsilon(h_n)^{-2}$ factor), the trace-operator-valued measure does not vanish and it must be dealt with using additional tools [Bur23]. Our result generalizes the estimate obtained by N. Burq and P. Gérard on the 2-dimensional torus in [BG20, Section 3A] (itself a generalization of [BZ15, Section 3]). The proof follows theirs.

Assume that $\|(h_n^2 \Delta + 1) u_n\|_{L^2} = o(h_n)$ and denote

$$(3.1) \qquad \epsilon(h_n) = \max\left( h_n^{\frac{1}{6}}, \left( \frac{\|(h_n^2 \Delta + 1) u_n\|}{h_n} \right)^{\frac{1}{6}} \right).$$

$\epsilon(h_n)$ satisfies:

$$(3.2) \qquad h_n^{-1} \epsilon^{-6}(h_n) \|(h_n^2 \Delta + 1) u_n\|_{L^2} \leq 1 \text{ and } \lim_{n \to +\infty} \epsilon(h_n) = 0.$$

The result we prove is the following:

**Proposition 3.1.** — *On the torus $\mathbb{T}^d$, assume that $\|u_n\|_{L^2(\mathbb{T}^d)} = \mathcal{O}(1)$ and (3.2) is satisfied, then there exists a constant $C > 0$ such that*

$$\forall n \in \mathbb{N}, \|u_n\|_{L^2(\{|x_j| \leq h_n^{\frac{1}{2}} \epsilon^{-2}(h_n)\})} \leq C \epsilon^{\frac{1}{2}}(h_n), \quad j = 1, \dots, d.$$

Note that $\epsilon(h_n) \geq h_n^{\frac{1}{6}}$ so the width of the slice tends to zero as $n$ goes to $+\infty$.

We prove the result for $j = 1$ without loss of generality to simplify notations. In this proof only, we denote $X = (x_1, y') \in \mathbb{T} \times \mathbb{T}^{d-1}$ and $\Xi = (\xi_1, \eta') \in \mathbb{R} \times \mathbb{R}^{d-1}$. Proposition 3.1 will result from the following statement:

**Proposition 3.2.** — *There exist some positive $C$ and $h_0$ such that for any $0 < h < h_0$, $1 \leq \beta \leq h^{-\frac{1}{2}}$ and any $(u, f) \in H^2 \times L^2$ satisfying*

$$(h^2(\partial_{x_1}^2 + \Delta_{y'}) + 1)u = f,$$

*the following estimate holds:*

$$(3.3) \quad \|u\|_{L^\infty(\{|x_1| \leq \beta h^{1/2}; L^2_{y'}\})}$$
$$\leq C \beta^{-\frac{1}{2}} h^{-\frac{1}{4}} \left( \|u\|_{L^2_{x_1,y'}(\{\beta h^{1/2} \leq |x_1| \leq 2\beta h^{1/2}\})} + h^{-1} \beta^2 \|f\|_{L^2_{x_1,y'}(\{|x_1| \leq 2\beta h^{1/2}\})} \right).$$

Let us show that Proposition 3.2 implies Proposition 3.1. We apply Hölder's inequality in the $x_1$ variable and Proposition 3.2 with $\beta = \epsilon^{-3}(h_n) \leq h_n^{-\frac{1}{2}}$:

$$\|u_n\|_{L^2(\{|x_1|\leq h_n^{1/2}\epsilon^{-2}(h_n)\}\})} \leq h_n^{\frac{1}{4}}\epsilon^{-1}(h_n)\|u_n\|_{L^\infty(\{|x_1|\leq h_n^{\frac{1}{2}}\epsilon^{-3}(h_n)\};L^2_{y'})}$$

$$\leq C\epsilon^{\frac{1}{2}}(h_n)\left(\|u_n\|_{L^2_{x_1,y'}(\{h_n^{1/2}\epsilon^{-3}(h_n)\leq|x_1|\leq 2h_n^{1/2}\epsilon^{-3}(h_n)\})}\right.$$

$$\left.+h_n^{-1}\epsilon^{-6}(h_n)\|f_n\|_{L^2_{x_1,y'}(\{|x_1|\leq 2h_n^{1/2}\epsilon^{-3}(h_n)\})}\right)$$

$$\leq C\epsilon^{\frac{1}{2}}(h_n)(\|u_n\|_{L^2}+h_n^{-1}\epsilon^{-6}(h_n)\|f_n\|_{L^2})$$

$$\leq 2C\epsilon^{\frac{1}{2}}(h_n).$$

We now prove Proposition 3.2. Denote $v$ (resp. $g$) the partial Fourier transform of $u$ (resp. $f$) with regards to the $y'$ variable. For a fixed $x_1$, the Plancherel equality gives

$$\|v(x_1,.)\|_{L^2_{\eta'}} = K\|u(x_1,.)\|_{L^2_{y'}},$$

where $K$ depends on the periods of the torus in the $y'$ directions. Inequality (3.3) is then equivalent to

$$(3.4) \quad \|v\|_{L^\infty(\{|x_1|\leq \beta h^{1/2}\};L^2_{\eta'})} \leq C\beta^{-\frac{1}{2}}h^{-\frac{1}{4}}\left(\|v\|_{L^2(\{\beta h^{1/2}\leq |x_1|\leq 2\beta h^{1/2}\};L^2_{\eta'})}\right.$$

$$\left.+h^{-1}\beta^2\|g\|_{L^2(\{|x_1|\leq 2\beta h^{1/2}\};L^2_{\eta'})}\right).$$

Besides, by the Minkowski inequality,

$$\|v\|_{L^\infty_{x_1};L^2_{\eta'}} \leq \|v\|_{L^2_{\eta'};L^\infty_{x_1}}$$

so it is enough to prove the following one-dimensional result:

**Proposition 3.3** ([BG20], Proposition 3.3). – *There exist some positive $C$ and $h_0$ such that for any $0 < h < h_0$, $\eta' \in \mathbb{R}$, $1 \leq \beta \leq h^{-\frac{1}{2}}$ and any $(v,g)$ satisfying*

$$\left(h^2\frac{d^2}{dx_1^2}+1-h^2|\eta'|^2\right)v = g,$$

*we have*

$$(3.5) \quad \|v\|_{L^\infty(\{|x_1|\leq \beta h^{\frac{1}{2}}\})} \leq C\beta^{-\frac{1}{2}}h^{-\frac{1}{4}}\left(\|v\|_{L^2(\{\beta h^{\frac{1}{2}}\leq|x|\leq 2\beta h^{\frac{1}{2}}\})}+h^{-1}\beta^2\|g\|_{L^2(\{|x|\leq 2\beta h^{\frac{1}{2}}\})}\right).$$

Integrating the square of (3.5) in the $\eta'$ variable gives (3.4).

The proof of the one-dimensional estimate (3.5) is exactly that of [BG20]. We include it for the sake of completeness. Performing the change of variables $x_1 = \beta h^{\frac{1}{2}}z$, it is equivalent to prove that any solutions of

$$\left(h\beta^{-2}\partial_z^2+1-h^2|\eta'|^2\right)v = g$$

satisfy

$$\|v\|_{L^\infty(\{|z|\leq 1\})} \leq C\left(\|v\|_{L^2(\{1\leq|z|\leq 2\})}+h^{-1}\beta^2\|g\|_{L^2(\{|z|\leq 2\})}\right).$$

Setting $\tau = \beta^2 h^{-1}(1-h^2|\eta'|^2)$, it suffices to prove the following lemma:

**Lemma 3.4** ([BG20], Lemma 3.4). — *There exists some positive $C$ such that for any $\tau \in \mathbb{R}$ and any solution $(v, k)$ over $(-2, 2)$ of*

$$(\partial_z^2 + \tau)v = k,$$

*there holds the inequality*

$$\|v\|_{L^\infty(-1,1)} \leq C\left(\|v\|_{L^2(\{1\leq|z|\leq 2\})} + \frac{1}{\sqrt{1+|\tau|}}\|k\|_{L^1(-2,2)}\right).$$

*Proof of the lemma.* Let $\chi \in C_c^\infty(-2, 2)$ be equal to 1 over $(-1, 1)$. $w = \chi v$ satisfies

(3.6) $$(\partial_z^2 + \tau)w = \chi k + 2\partial_z(\chi' v) - \chi'' v.$$

We distinguish two regimes:

Elliptic regime: $\tau \leq -1$. We multiply (3.6) by $w$ and integrate by parts to get

$$\|\partial_z w\|_{L^2(-2,2)}^2 + |\tau|\|w\|_{L^2(-2,2)}^2 = -(\chi k - \chi'' v, w)_{L^2} + 2(\chi' v, \partial_z w)_{L^2}.$$

Thus,

$$\|\partial_z w\|_{L^2(-2,2)}^2 + |\tau|\|w\|_{L^2(-2,2)}^2 \leq$$
$$C\left(\|k\|_{L^1(-2,2)}\|w\|_{L^\infty} + \|v\|_{L^2(\{1\leq|z|\leq 2\})}(\|w\|_{L^2(\{1\leq|z|\leq 2\})} + \|\partial_z w\|_{L^2(-2,2)})\right).$$

By the Gagliardo-Nirenberg inequality in dimension 1,

$$\|w\|_{L^\infty} \leq C\|\partial_z w\|_{L^2}^{\frac{1}{2}}\|w\|_{L^2}^{\frac{1}{2}} = C|\tau|^{-\frac{1}{4}}\|\partial_z w\|_{L^2}^{\frac{1}{2}}|\tau|^{\frac{1}{4}}\|w\|_{L^2}^{\frac{1}{2}}.$$

Set $A = \|\partial_z w\|_{L^2(-2,2)} + |\tau|^{\frac{1}{2}}\|w\|_{L^2(-2,2)}$, then applying the Gagliardo-Nirenberg inequality again gives

$$A^2 \leq C\left(\|k\|_{L^1(-2,2)}|\tau|^{-\frac{1}{4}}A + \|v\|_{L^2(\{1\leq|z|\leq 2\})}A\right).$$

Applying it a third time yields

$$\|v\|_{L^\infty(-1,1)} \leq \|w\|_{L^\infty} \leq C|\tau|^{-\frac{1}{4}}A \leq C\left(\|v\|_{L^2(\{1\leq|z|\leq 2\})} + |\tau|^{-\frac{1}{2}}\|k\|_{L^1(-2,2)}\right).$$

Since $|\tau| \geq 1$, we get the desired inequality.

Hyperbolic regime: $\tau \geq -1$. Set $\sigma = \sqrt{\tau} \in \mathbb{R}^+ \cup i[0, 1]$. Integrating (3.6) gives

$$w(x) = \int_{z=-2}^x g(z) \int_{y=z}^x e^{i\sigma(2y-x-z)} dy dz,$$

where $g = \chi k - \chi'' v + 2\partial_z(\chi' v) = g_1 + \partial_z g_2$. Since for any $x, z \in [-2, 2]$

$$\left|\int_{y=z}^x e^{i\sigma(2y-x-z)} dy\right| \leq \frac{C}{1+|\sigma|} = \frac{C}{1+|\tau|^{\frac{1}{2}}},$$

the contribution of $g_1$ is bounded uniformly by

$$\frac{C}{1+|\tau|^{\frac{1}{2}}}\left(\|\chi k\|_{L^1(-2,2)} + \|v\|_{L^1(\{1\leq|z|\leq 2\})}\right).$$

Integrating by parts, we get that the contribution of $\partial_z g_2$ is bounded by

$$C\|\chi' v\|_{L^1(-2,2)} \leq C'\|v\|_{L^1(\{1\leq|z|\leq 2\})}.$$

Summing both contributions gives the result. This concludes the proof of Proposition 3.1. □

# 4. Proof of Theorem 4

This section is dedicated to the proof of Theorem 4, namely that stabilization occurs when every geodesic either intersects the interior of the damped zone or is damped in every normal direction. The first subsection is dedicated to the proof of $L^2$-boundedness of pseudodifferential operators and a Gårding inequality for our 2-microlocal calculus. Well-known semiclassical analysis techniques are used to give an appropriate and precise framework to define 2-microlocal measures. In the second subsection, we study the 2-microlocal measure and prove uniform stabilization for a model damping. In the third, we reduce the general case to this model. The last subsection contains the proof of a geometric lemma used in the second one.

## 4.1. Second microlocal calculus.

We start by constructing the 2-microlocal pseudodifferential calculus. We deal with more general symbols than in [BG20], as we only assume decay in the $\zeta$ variable rather than polyhomogeneity. This weakened assumption will be used in the third microlocalization in Section 5.

We introduce symbol classes and the second quantization for these symbols, then show that this quantization provides a pseudodifferential calculus in $\epsilon^2(h)$ with sufficient tools to construct 2-microlocal measures. Our symbols take variables in $\mathbb{R}^d \times \mathbb{R}^d \times \mathbb{R}^q \times \mathbb{R}^q$. Case $q = d - 1$ will be used in the present section, and $d = 3$, $q \in \{1, 2\}$ in Section 5.

**Definition 4.1** (Symbol classes). — (1) We define $S^0$ to be the class of smooth functions $b$ of $(X, \Xi, z, \zeta) \in \mathbb{R}^d \times \mathbb{R}^d \times \mathbb{R}^q \times \mathbb{R}^q$ which are compactly supported in the $(X, \Xi)$ variables and satisfy the following decay estimate in the $\zeta$ variable:

$$\forall \alpha, \beta \in \mathbb{N}^q, \quad \exists C_{\alpha,\beta} \quad s.\ t. \quad \sup_{(X,\Xi,z,\zeta) \in \mathbb{R}^{2d} \times \mathbb{R}^{2q}} |\partial_z^\alpha \partial_\zeta^\beta b(X, \Xi, z, \zeta)| \leq C_{\alpha,\beta} \langle \zeta \rangle^{-|\beta|}$$

where $\langle . \rangle$ denotes the usual $\langle \zeta \rangle = (1 + |\zeta|^2)^{\frac{1}{2}}$.

(2) For any integer $m \geq 0$, we define $S_H^m$ as the set of smooth functions which are compactly supported in the $(X, \Xi)$ variables and polyhomogeneous of degree $m$ with respect to the $(z, \zeta)$ variables with limits in the radial direction:

$$\lim_{r \to +\infty} \frac{1}{r^m} a\left(X, \Xi, \frac{(rz, r\zeta)}{\|(z,\zeta)\|}\right) = \tilde{a}\left(X, \Xi, \frac{(z,\zeta)}{\|(z,\zeta)\|}\right).$$

We have $S_H^0 \subset S^0$. Besides, when $m = 0$ functions in $S_H^0$ are identified with smooth compactly supported functions on $\mathbb{R}^{2d} \times \overline{B(0,1)}_{\tilde{z},\tilde{\zeta}}$ via the change of variables

$$(z, \zeta) \mapsto (\tilde{z}, \tilde{\zeta}) = \frac{(z,\zeta)}{\sqrt{1+|z|^2+|\zeta|^2}}.$$

Consider some function $\epsilon(h)$ satisfying

(4.1) $$\lim_{h \to 0} \epsilon(h) = 0, \quad \epsilon(h) \geq h^{\frac{1}{2}}, \quad h^{\frac{1}{2}}\epsilon^{-2}(h) \to_{h \to 0} 0.$$

Given a symbol $b$ belonging to $S_H^m$ or $S^0$, its second quantization is defined by:

(4.2) $$\mathrm{Op}_h^{(\epsilon)}(b)u(X) = \frac{1}{(2\pi)^d} \int_{\mathbb{R}_Y^d \times \mathbb{R}_\Xi^d} e^{i(X-Y)\cdot\Xi} b\left(X, h\Xi, \frac{\epsilon(h)}{h^{\frac{1}{2}}}x, \epsilon(h)h^{\frac{1}{2}}\xi\right) u(Y) dY d\Xi$$

where $X = (x, x'), \Xi = (\xi, \xi')$ belong to $\mathbb{R}^q \times \mathbb{R}^{d-q}$. In other words,

$$\text{Op}_h^{(\epsilon)}(b) = b\left(x, x', hD_x, hD_{x'}, \frac{\epsilon(h)}{h^{\frac{1}{2}}}x, \epsilon(h)h^{\frac{1}{2}}D_x\right). \tag{4.3}$$

We first prove uniform boundedness over $L^2(\mathbb{R}^d)$ and over the space $L^2_{ul}$ defined in (4.4). Boundedness over $L^2_{ul}$ will allow to define 2-microlocal measures for periodic functions over $\mathbb{R}^d$.

**Proposition 4.2.** – *Consider some symbol $b \in S^0$, then there exists some positive $h_0$ such that the following statements hold:*

(1) $\text{Op}_h^{(\epsilon)}(b)$ *is a bounded operator over $L^2(\mathbb{R}^d)$ uniformly wrt $h \in (0, h_0]$.*

(2) *We introduce a partition of unity $\chi \in C^\infty(\mathbb{R}^d)$ such that $\text{supp}(\chi) \subset [-1,1]^d$, $0 \leq \chi \leq 1$ and $\sum_{p \in \mathbb{Z}^d} \chi_p = 1$ where $\chi_p = \chi(\cdot - p)$. We define the set of uniformly locally $L^2$ functions over $\mathbb{R}^d$ by*

$$L^2_{ul} = \left\{ u \in L^2_{loc}, \ \sup_{p \in \mathbb{Z}^d} \|\chi_p u\|_{L^2(\mathbb{R}^d)} < +\infty \right\} \tag{4.4}$$

*endowed with the norm $\|u\|_{L^2_{ul}} = \sup_{p \in \mathbb{Z}^d} \|\chi_p u\|_{L^2(\mathbb{R}^d)}$. Then $\text{Op}_h^{(\epsilon)}(b)$ is a bounded operator over $L^2_{ul}(\mathbb{R}^d)$ uniformly wrt $h \in (0, h_0]$.*

*Proof.* We first prove uniform boundedness over $L^2$. Setting $\tilde{x} = \frac{\epsilon(h)}{h^{\frac{1}{2}}}x$, $\tilde{y} = \frac{\epsilon(h)}{h^{\frac{1}{2}}}y$, and $\tilde{\xi} = \epsilon(h)h^{\frac{1}{2}}\xi$ in (4.2), we obtain that

$$\text{Op}_h^{(\epsilon)}(b) = T_h^* b\left(\frac{h^{\frac{1}{2}}}{\epsilon(h)}x, x', \epsilon(h)h^{\frac{1}{2}}D_x, hD_{x'}, x, \epsilon^2(h)D_x\right) T_h \tag{4.5}$$

where $T_h$ is the unitary operator defined by $T_h u(x, x') = \left(\frac{h^{\frac{1}{2}}}{\epsilon(h)}\right)^{\frac{q}{2}} u\left(\frac{h^{\frac{1}{2}}}{\epsilon(h)}x, x'\right)$. Since $\epsilon(h) \geq h^{\frac{1}{2}}$, the $h$-dependent symbol

$$\check{b}(X, \Xi) = b\left(\frac{h^{\frac{1}{2}}}{\epsilon(h)}x, x', \epsilon(h)h^{\frac{1}{2}}\xi, h\xi', x, \epsilon^2(h)\xi\right)$$

has bounded derivatives uniformly wrt $h$ small enough. By the Calderón-Vaillancourt theorem [Zwo12, Theorem 5.1], $\text{Op}_1(\check{b})$ is bounded over $L^2$, which proves uniform $L^2$-boundedness for $\text{Op}_h^{(\epsilon)}(b)$ by unitary conjugation.

To generalize the result to $L^2_{ul}$, consider the symbol $\check{b}$ defined above and write, up to unitary conjugation:

$$\chi_r \text{Op}_h^{(\epsilon)}(b) u(X) = \chi_r(X) \sum_{p,q \in \mathbb{Z}^d} \frac{1}{(2\pi)^d} \int_{\mathbb{T}^d \times \mathbb{R}^d} e^{i(X-Y).\Xi} \chi_p(X) \check{b}(X, \Xi) \chi_q(Y) u(Y) dY d\Xi.$$

Our goal is to obtain some bound on $\|\chi_r \text{Op}_h^{(\epsilon)}(b) u\|_{L^2}$ that is uniform wrt $h \in (0, h_0]$ and $r \in \mathbb{Z}^d$. If $|p - r| \geq 2$ then $\text{supp}(\chi_p) \cap \text{supp}(\chi_r)$ has measure zero, hence the the sum over $p$ is finite. Then, whenever $|p - q| \geq 2$, we integrate by parts in the $\Xi$ variable using that $\check{b}$ is compactly supported and that for any integer $N$

$$\left(\frac{1 - \Delta_\Xi}{1 + |X - Y|^2}\right)^N \left(e^{i(X-Y).\Xi}\right) = e^{i(X-Y).\Xi}.$$

Since $X$ is localized near $p$ and $Y$ is localized near $q$, when $|p - q| \geq 2$ there exists some $c > 0$ such that
$$\frac{1}{c}\frac{1}{1+|p-q|} \leq \frac{1}{1+|X-Y|} \leq \frac{c}{1+|p-q|}.$$
The sum over $q$ is thus convergent for every fixed $p$ and $N$ large enough. The Calderón-Vaillancourt theorem then allows to conclude, as

$$\|\chi_r \operatorname{Op}_h^{(\epsilon)}(b)u\| \leq C \sum_{|p-r|\leq 1} \sum_{q \in \mathbb{Z}^d} \left\| \frac{\chi_p}{(1+|p-q|)^N} \operatorname{Op}_h^{(\epsilon)}((1-\Delta_\Xi)^N b) \chi_q u \right\|_{L^2}$$

$$\leq C \sum_{|p-r|\leq 1} \left( \sum_{|p-q|\geq 2} \frac{1}{(1+|p-q|)^N} \|\chi_q u\|_{L^2} + \sum_{|p-q|<2} \|\chi_q u\|_{L^2} \right)$$

$$\leq C \left( \sum_{q' \in \mathbb{Z}^d} \frac{1}{(1+|q'|)^N} \right) \|u\|_{L^2_{ul}}.$$

$\square$

We now prove Gårding's weak inequality with small parameter $\epsilon^2(h)$. The result is the following:

**Proposition 4.3.** — *Let $b \in S^0$ be non-negative. Then, for any $\alpha > 0$, there exist some $h_0(\alpha) > 0$ and $C(\alpha) > 0$ such that for any $h \in (0, h_0]$, $u \in L^2(\mathbb{R}^d)$:*

(4.6)
$$\operatorname{Re}(b\left(\frac{h^{\frac{1}{2}}}{\epsilon(h)}x, x', \epsilon(h)h^{\frac{1}{2}}D_x, hD_{x'}, x, \epsilon^2(h)D_x\right)u, u)_{L^2} \geq -(\alpha + C\epsilon^2(h))\|u\|_{L^2}^2.$$
$$\left| \operatorname{Im}(b\left(\frac{h^{\frac{1}{2}}}{\epsilon(h)}x, x', \epsilon(h)h^{\frac{1}{2}}D_x, hD_{x'}, x, \epsilon^2(h)D_x\right)u, u)_{L^2} \right| \leq C\epsilon^2(h)\|u\|_{L^2}^2.$$

*By unitary conjugation, we conclude that the same inequality holds for the operator*
$$b\left(x, x', hD_x, hD_{x'}, \frac{\epsilon(h)}{h^{\frac{1}{2}}}x, \epsilon(h)h^{\frac{1}{2}}D_x\right).$$

*Proof.* We first note that $b\left(\frac{h^{\frac{1}{2}}}{\epsilon(h)}x, x', \epsilon(h)h^{\frac{1}{2}}D_x, hD_{x'}, x, \epsilon^2(h)D_x\right)$ is the $\epsilon^2(h)$-quantization of the non-negative $h$-dependent symbol

(4.7)
$$\tilde{b}_h(X, \Xi) = b\left(\frac{h^{\frac{1}{2}}}{\epsilon(h)}x, x', \frac{h^{\frac{1}{2}}}{\epsilon(h)}\xi, \frac{h}{\epsilon^2(h)}\xi', x, \xi\right).$$

Fix some $\alpha > 0$ and consider the symbol $c_h = \sqrt{\tilde{b}_h + \alpha}$. Since $h^{\frac{1}{2}}\epsilon^{-1}(h) \leq 1$, all the derivatives of $\tilde{b}_h$ are bounded over $\mathbb{R}^d \times \mathbb{R}^d$ uniformly wrt $h$ small enough and derivatives of order $m$ in the $\Xi$ variable decay at rate $\langle \Xi \rangle^{-m}$. Thus, $c_h$ also has bounded derivatives and a similar rate of decay. Since $c_h$ is real-valued, we get that
$$\operatorname{Op}_{\epsilon^2(h)}(c_h)^* = \operatorname{Op}_{\epsilon^2(h)}(c_h) + \mathcal{O}_{\mathcal{L}(L^2)}\left(\epsilon^2(h)\right).$$
We then apply the usual composition formulas (see for example [Zwo12, Theorem 4.14] which generalizes to our framework) to get
$$(\operatorname{Op}_{\epsilon^2(h)}(\tilde{b}_h)u, u)_{L^2} + \alpha\|u\|_{L^2}^2 = \|\operatorname{Op}_{\epsilon^2(h)}(c_h)u\|_{L^2}^2 + \mathcal{O}(\epsilon^2(h))\|u\|_{L^2}^2.$$

Taking the real and imaginary parts then gives the desired inequalities (4.6) for $h_0(\alpha)$ small enough and some constant $C(h_0, \alpha)$. □

4.2. **The model case: 4 damped prisms in a tunnel.** Using the 2-microlocal calculus for symbols in $S_H^0$, we turn to the study of the concentration properties of a sequence $(u_n)$ of quasimodes over the three-dimensional torus $\mathbb{T}^3$, under a specific choice of damping function $a$. We will then generalize this model case to higher dimensions and reduce the general case to it.

We identify $\mathbb{T}^3$ with $[-1, 1]^3$ where opposite faces are identified. The damping zone $\{a = 1\}$ is represented in Figure 4.1. More explicitly, denoting $X = (x_1, x_2, y)$ and $\Xi = (\xi_1, \xi_2, \eta)$, $a$ is the characteristic function of the set

$$
\begin{aligned}
&\{|x_1| \geq \frac{1}{2} \text{ or } |x_2| \geq \frac{1}{2}\} \\
&\cup \{-1 < y < -\frac{1}{2}, x_1 > 0, x_2 < \alpha_R x_1\} \cup \{-\frac{1}{2} < y < 0, x_2 > 0, x_1 > -\alpha_T x_2\} \\
&\cup \{0 < y < \frac{1}{2}, x_1 < 0, x_2 > \alpha_L x_1\} \cup \{\frac{1}{2} < y < 1, x_2 < 0, x_1 < -\alpha_B x_2\}.
\end{aligned}
$$
(4.8)

The first set denotes a neighborhood of the lateral faces of the cube as represented in the central picture of Figure 4.1 and the four other sets correspond to the prisms represented on the right, from the back to the front. $\alpha_L, \alpha_R, \alpha_T, \alpha_B > 0$ are four fixed positive coefficients which should be thought of as small. The case $\alpha_L = \alpha_R = \alpha_T = \alpha_B = 0$ (when the triangle-shaped damped sections vanish and the four prisms are reduced to cubes) is dealt with in section 5.

In this situation, the only two bicharacteristics that do not intersect the interior of $\text{supp}(a)$ are

$$\{x_1 = x_2 = 0, \xi_1 = \xi_2 = 0, \eta = \pm 1\},$$

which are represented in red in the figure (as one geodesic, traveled in two directions). Thus, the first microlocal measure $\nu$ is supported in the union of these two bicharacteristics. The corresponding geodesic satisfies assumption (1.3) if all $\alpha_{L,R,T,B}$ coefficients are positive but violates it if any of the coefficients vanishes.

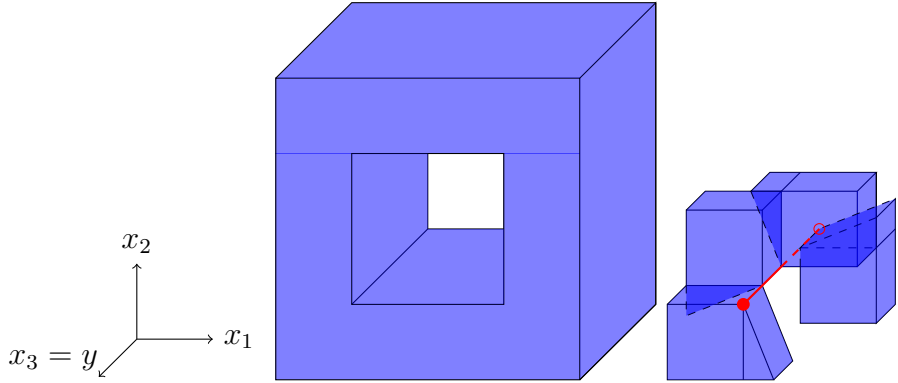

FIGURE 4.1. Illustration of set (4.8). The coordinate system is represented on the left. The dual coordinates are denoted $(\xi_1, \xi_2, \xi_3 = \eta)$. We also denote $x = (x_1, x_2)$, $\xi = (\xi_1, \xi_2)$. In the center is the torus $\mathbb{T}^3$ identified with the unit cube $[-1, 1]^3$. We color in blue the zone $\{|x_1| \geq \frac{1}{2}\} \cup \{|x_2| \geq \frac{1}{2}\}$ where $a = 1$, thus damping the lateral faces. This leaves an undamped tunnel at the center of the cube. The picture on the right fits inside this tunnel but is represented separately for better legibility. On the right, we represent in blue the set where $a = 1$ inside $\{|x_1| \leq \frac{1}{2}, |x_2| \leq \frac{1}{2}\}$. The red geodesic $\{x = 0\}$ is the only one which does not encounter the interior of $\text{supp}(a)$.

We now perform the second microlocalization around this geodesic. Consider a sequence $(u_n)$ of functions over $\mathbb{T}^3$ satisfying

(4.9) $$(h_n^2 \Delta + 1)u_n = \mathcal{O}(1)_{L^2}.$$

We identify $u_n \in L^2(\mathbb{T}^3)$ with its periodic extension over $\mathbb{R}^3$. In particular, the periodic extension of $u_n$ belongs to $L^2_{ul}$ as defined in (4.4). Using the $L^2_{ul}$-boundedness of operators defined by (4.2) and Gårding's inequality, we can extract a subsequence (also denoted $(u_n)$) such that there exists a positive measure $\tilde{\mu}$ on $T^*\mathbb{R}^3 \times \overline{N}$ satisfying for any polyhomogeneous symbol $b \in S_H^0$,

$$\lim_{n \to +\infty} (\text{Op}_{h_n}^{(\epsilon)}(b)u_n, u_n)_{L^2} = \langle \tilde{\mu}, \tilde{b} \rangle.$$

Here, $\overline{N}$ denotes the sphere compactification of $\mathbb{R}^4_{z,\zeta}$ and $\tilde{b}$ is a continuous function on $T^*\mathbb{R}^3 \times \overline{N}$ which is defined by the value of $b$ in the interior and by

$$\tilde{b}(x, y, \xi, \eta, \tilde{z}, \tilde{\zeta}) = \lim_{r \to +\infty} b(x, y, \xi, \eta, r\tilde{z}, r\tilde{\zeta})$$

on the sphere at infinity. Since the $u_n$ are periodic, the measure $\tilde{\mu}$ is also periodic in the $X$ variable, hence it naturally defines a measure $\mu$ on $T^*\mathbb{T}^3 \times \overline{N}$.

The rest of the subsection is dedicated to studying such a 2-microlocal measure $\mu$ for the sequence $(u_n)$ introduced in Section 2. Recall that $(u_n)$ satisfies $\|u_n\|_{L^2(\mathbb{T}^d)} = 1$ and the asymptotics:

$$(h_n^2 \Delta + 1)u_n = \mathcal{O}(h_n \epsilon(h_n))_{L^2} \quad ; \quad au_n = o(1)_{L^2}.$$

The properties of the second microlocal measure $\mu$ are gathered in the following statement :

**Proposition 4.4.** –   (1) *Assume only that*

$$(h_n^2\Delta + 1)u_n = \mathcal{O}(1)_{L^2},$$

*then $\mu(T^*\mathbb{T}^3 \times \overline{N}) = 1$.*

(2) *The first microlocal measure $\nu$ is the projection of $\mu$ onto the $(x, y, \xi, \eta)$ variables. Assume besides that*

$$(h_n^2\Delta + 1)u_n = o(h_n)_{L^2}, \quad au_n = o(1)_{L^2},$$

*then $\mu$ is supported 1-microlocally in the red bicharacteristics of Figure 4.1,*

$$\{x = 0, \quad \xi = 0, \quad \eta = \pm 1\}.$$

(3) *Assume now that*

$$(h_n^2\Delta + 1)u_n = \mathcal{O}(h_n\epsilon(h_n))_{L^2},$$

*then $\mu$ is supported in the sphere at infinity in the $(z, \zeta)$ variables.*

(4) *Let $R$ be one of the four prisms*

$$\{-1 < y < -\frac{1}{2}, x_1 > 0, x_2 < \alpha_R x_1\}, \{-\frac{1}{2} < y < 0, x_2 > 0, x_1 > -\alpha_T x_2\},$$

$$\{0 < y < \frac{1}{2}, x_1 < 0, x_2 > \alpha_L x_1\}, \{\frac{1}{2} < y < 1, x_2 < 0, x_1 < -\alpha_B x_2\}$$

*(see Figure 4.1), then at every point of $\partial R$ $\mu$ vanishes 2-microlocally in the direction of the polyhedron $R$. Precisely (from back to front along the red geodesic of Figure 4.1):*

(4.10)
$$\mu(\{x = 0, \quad \xi = 0, \quad -1 < y < -\frac{1}{2}, \quad \eta = \pm 1, \quad z_1 > 0, \quad z_2 < \alpha_R z_1\}) = 0,$$
$$\mu(\{x = 0, \quad \xi = 0, \quad -\frac{1}{2} < y < 0, \quad \eta = \pm 1, \quad z_2 > 0, \quad z_1 > -\alpha_T z_2\}) = 0,$$
$$\mu(\{x = 0, \quad \xi = 0, \quad 0 < y < \frac{1}{2}, \quad \eta = \pm 1, \quad z_1 < 0, \quad z_2 > \alpha_L z_1\}) = 0,$$
$$\mu(\{x = 0, \quad \xi = 0, \quad \frac{1}{2} < y < 1, \quad \eta = \pm 1, \quad z_2 < 0, \quad z_1 < -\alpha_B z_2\}) = 0.$$

(5) *The measure $\mu$ satisfies the conservation law:*

(4.11) $$(\eta\partial_y + \zeta.\partial_z)\mu = 0$$

*where $\zeta.\partial_z$ is a vector field on the sphere at infinity of $\mathbb{R}^4_{z,\zeta}$.*

*Proof.* The first point is a direct consequence of the $h_n$-oscillation property proved in section 2.

To prove the second statement, consider a symbol $b(X, \Xi)$ which does not depend on the 2-microlocal variables $z, \zeta$, then $\langle \nu, b\rangle = \langle \mu, b \circ \pi\rangle$ by definition of both measures. Said otherwise, $\nu = \pi_*\mu$ where

$$\pi : (X, \Xi, z, \zeta) \mapsto (X, \Xi).$$

The rest of the second point comes from Corollary 2.2.

To prove the third point, we take some function $\chi \in C_c^\infty(T^*\mathbb{T}^3 \times \mathbb{R}^4_{z,\zeta})$. $\chi$ belongs to $S_H^0$ and there exists some positive $A$ such that $\mathrm{supp}(\chi) \subset \{|z| \leq A\}$. Thus,

$$
\begin{aligned}
(4.12) \quad |(\mathrm{Op}_{h_n}^{(\epsilon)}(\chi)u_n, u_n)_{L^2}| &= |(\mathrm{Op}_{h_n}^{(\epsilon)}(\chi)u_n, \mathbb{1}_{\{|x| \leq Ah_n^{\frac{1}{2}}\epsilon^{-1}(h_n)\}} u_n)_{L^2}| \\
&\leq C\|u_n\|_{L^2}\|u_n\|_{L^2(\{|x| \leq h_n^{\frac{1}{2}}\epsilon^{-2}(h_n)\})} \\
&\leq C\|u_n\|_{L^2(\{|x_1| \leq h_n^{\frac{1}{2}}\epsilon^{-2}(h_n)\})} \\
&\leq C\epsilon^{\frac{1}{2}}(h_n) \to_{n \to +\infty} 0.
\end{aligned}
$$

where we have applied the $L^2$-boundedness of $\mathrm{Op}_{h_n}^{(\epsilon)}(\chi)$ and Proposition 3.1. Thus,

$$\langle \mu, \chi \rangle = 0,$$

which proves that $\mu$ is supported in the sphere at infinity in the $(z, \zeta)$ variables. Note that it is enough for the non-concentration estimate to hold in the cylinder $\{|x| \leq h_n^{\frac{1}{2}}\epsilon^{-2}(h_n)\}$ centered around the geodesic rather than the slice $\{|x_1| \leq h_n^{\frac{1}{2}}\epsilon^{-2}(h_n)\}$. This fact will be used in the general case at the end of the next subsection.

We now come to show the microlocal vanishings. We show that $\mu$ vanishes over

$$\{x = 0, \xi = 0, \frac{1}{2} < y < 1, \eta = +1, z_2 < 0, z_1 < -\alpha_B z_2\}.$$

Fix some small enough $\delta_0 > 0$ and define the following cutoffs:

- $\psi \in C^\infty(\mathbb{R})$, equal to 0 over $(-\infty, 1]$ and 1 over $[2, +\infty)$,
- $\chi \in C_c^\infty(-1, 1)$ equal to 1 over $(-\frac{1}{2}, \frac{1}{2})$, and
- $\tilde{\chi} \in C_c^\infty(\frac{1}{2}, 1)$ equal to 1 over $(\frac{1}{2} + \delta_0, 1 - \delta_0)$.

Consider the symbol $b \in S_H^0$ defined by

$$(4.13) \quad b(x, y, \xi, \eta, z, \zeta) =$$
$$\chi\left(\frac{2|x|}{\delta_0}\right) \tilde{\chi}(y) \chi(|\xi|) \chi(\eta - 1) \psi\left(-\frac{z_2}{\delta_0 |\zeta|}\right) \psi\left(\frac{-z_1 - \alpha_B z_2}{\delta_0 |\zeta|}\right) \psi\left(\frac{|z|^2 + |\zeta|^2}{\delta_0}\right).$$

Cutoff $\chi\left(\frac{2|x|}{\delta_0}\right)\tilde{\chi}(y)$ is supported in $B(0, \frac{\delta_0}{2})_x \times (\frac{1}{2}, 1)_y$ and cutoff $\psi\left(-\frac{z_2}{\delta_0|\zeta|}\right)\psi\left(\frac{-z_1 - \alpha_B z_2}{\delta_0 |\zeta|}\right)$ localizes in the $\{x_2 < 0, x_1 < -\alpha_B x_2\}$ zone, thus

$$\mathbb{1}_{\{x \in B(0, \frac{\delta_0}{2}), x_2 < 0, x_1 < -\alpha_B x_2\}} \mathbb{1}_{\{y \in (\frac{1}{2}, 1)\}} \mathrm{Op}_{h_n}^{(\epsilon)}(b) = \mathrm{Op}_{h_n}^{(\epsilon)}(b),$$

wherefrom

$$
\begin{aligned}
(4.14) \quad |\langle \mathrm{Op}_{h_n}^{(\epsilon)}(b)u_n, u_n \rangle_{L^2}| &= |\langle \mathrm{Op}_{h_n}^{(\epsilon)}(b)u_n, \mathbb{1}_{\{x \in B(0, \frac{\delta_0}{2}), x_2 < 0, x_1 < -\alpha_B x_2\}}\mathbb{1}_{\{y \in (\frac{1}{2}, 1)\}} u_n \rangle_{L^2}| \\
&= |\langle \mathrm{Op}_{h_n}^{(\epsilon)}(b)u_n, \mathbb{1}_{\{x \in B(0, \frac{\delta_0}{2}), x_2 < 0, x_1 < -\alpha_B x_2\}}\mathbb{1}_{\{y \in (\frac{1}{2}, 1)\}} au_n \rangle_{L^2}| \\
&\leq C\|u_n\|_{L^2}\|au_n\|_{L^2} \to_{n \to +\infty} 0.
\end{aligned}
$$

Thus $\langle \mu, b \rangle = 0$. $\mu$ being a positive measure, this implies $\mu(\{b = 1\}) = 0$ and

$$\mu\left(\left\{x = 0, \xi = 0, y \in \left[\frac{1}{2} + \delta_0, 1 - \delta_0\right], \eta = 1, z_2 < -2\delta_0|\zeta|, -(z_1 + \alpha_B z_2) > 2\delta_0|\zeta|\right\}\right) = 0.$$

Taking the limit $\delta_0 \to 0$ then yields the result.

The other 2-microlocal vanishings are obtained by changing the cutoffs in the $y$, $\eta$ and $\frac{z}{|\zeta|}$ variables in (4.13).

To prove the last statement, we compute the bracket $[h_n^2 \Delta + 1, \text{Op}_{h_n}^{(\epsilon)}(q)]$ for some $q \in S_H^0$ using the second quantization formula (4.2). We get

$$[h\Delta, \text{Op}_h^{(\epsilon)}(q)] = \text{Op}_h^{(\epsilon)}\left((2i\Xi.\partial_X + 2i\zeta.\partial_z + h\Delta_X + 2h\epsilon(h)h^{-\frac{1}{2}}(\partial_x.\partial_z) + h(\epsilon(h)h^{-\frac{1}{2}})\Delta_z)q\right)$$

so that

$$(4.15) \quad \frac{1}{2ih_n}[h_n^2\Delta + 1, \text{Op}_{h_n}^{(\epsilon)}(q)] = \frac{1}{2i}[h_n\Delta, \text{Op}_{h_n}^{(\epsilon)}(q)]$$
$$= \text{Op}_{h_n}^{(\epsilon)}((\xi.\partial_x + \eta\partial_y + \zeta.\partial_z)q) - i\frac{h_n}{2}\text{Op}_{h_n}^{(\epsilon)}(\Delta_X q) - i\frac{h_n}{2}(\epsilon(h_n)h_n^{-\frac{1}{2}})\text{Op}_{h_n}^{(\epsilon)}((\partial_x.\partial_z)q)$$
$$- i\frac{h_n}{2}(\epsilon(h_n)h_n^{-\frac{1}{2}})^2\text{Op}_{h_n}^{(\epsilon)}(\partial_z^2 q).$$

Unfolding the bracket and using that $(h_n^2\Delta + 1)u_n = o(h_n)$, we obtain

$$\frac{1}{2ih_n}([h_n^2\Delta + 1, \text{Op}_{h_n}^{(\epsilon)}(q)]u_n, u_n)_{L^2} \to_{n\to+\infty} 0,$$

so

$$(4.16) \quad (\text{Op}_{h_n}^{(\epsilon)}((\xi.\partial_x + \eta\partial_y + \zeta.\partial_z)q)u_n, u_n)_{L^2} = o(1).$$

Since the measure $\mu$ is supported in $\{\xi = 0\}$, conservation law (4.11) follows. □

We can then conclude the contradiction argument by studying the dynamics of the vector field $\zeta.\partial_z$ over the sphere at infinity. The corresponding result is illustrated in Figure 4.3 of subsection 4.4 and summed up in the following lemma. We postpone the proof of the lemma to that subsection.

**Lemma 4.5.** – *Consider an initial point $(z_0, \zeta_0)$, $|z_0|^2 + |\zeta_0|^2 = 1$, in the sphere at infinity $\mathbb{S}_\infty^{2d-3}$ of $\mathbb{R}^{2d-2}$, and denote $\phi^s(z_0, \zeta_0)$ the flow of $\zeta.\partial_z$ at time $s$ starting from this initial point, then*

$$(4.17) \quad \phi^s(z_0, \zeta_0) = \frac{1}{(|z_0 + s\zeta_0|^2 + |\zeta_0|^2)^{\frac{1}{2}}}(z_0 + s\zeta_0, \zeta_0).$$

*In other words, the flow of of $\zeta.\partial_z$ over the sphere at infinity is the projection of the flow at finite distance from the origin, as depicted in Figure 4.3.*

*In particular, if $\zeta_0 \neq 0$ then*

$$(4.18) \quad \phi^s(z_0, \zeta_0) \to_{s\to+\infty} \left(\frac{\zeta_0}{|\zeta_0|}, 0\right),$$

*and $(z_0, \zeta_0)$ is a fixed point of the flow iff $\zeta_0 = 0$.*

Now, by point (1) of Proposition 4.4, $\mu$ has non-empty support. Consider a point

$$(x = 0, y, \xi = 0, \eta = \pm 1, z, \zeta) \in \text{supp}(\mu),$$

with $(z, \zeta)$ belonging to the sphere at infinity of $\mathbb{R}^4$. By the conservation law (4.11), the point

$$(x = 0, y + s\eta, \xi = 0, \eta = \pm 1, \phi^s(z, \zeta))$$

also belongs to supp($\mu$). The flow $\phi_s(z,\zeta)$ converges to some $(z_\infty, 0)$ with $z_\infty \neq 0$, and by the fourth point $(z_\infty, 0)$ belongs to at least one open set of the sphere at infinity where $\mu$ vanishes 2-microlocally along the geodesic. Thus, $\mu$ has support in an open set where it vanishes, which gives a contradiction. This concludes the proof of stabilization for the model damping in Figure 4.1.

4.2.1. *Generalization to higher dimensions.* Since both the non-concentration estimate of Proposition 3.1 and Lemma 4.5 hold for any $d \geq 3$, this model case can be generalized to higher dimensions. Assume that $d \geq 3$. We consider a damping $a$ over the $d-$dimensional torus $\mathbb{T}^d \simeq [-1,1]^d$ such that:

- $a = 1$ outside of a tunnel, namely in the zone $\cup_{i=1}^{d-1}\{|x_i| \geq \frac{1}{2}\}$. By Corollary 2.2, the first microlocal measure is then supported in the one-directional set of bicharacteristics $\{\xi_i = 0, 1 \leq i \leq d-1, \xi_d = \pm 1\}$.
- Inside that tunnel, $a = 1$ over a finite union of polyhedrons such that the geodesic $\{x_1 = \cdots = x_{d-1} = 0\}$ satisfies Assumption (1.3) and is the only geodesic not entering Int(supp($a$)). For example, one can choose $a$ to be equal to 1 over the set

$$\bigcup_{i=1}^{d-1}\left\{x_i < 0, \quad -1+\frac{i-1}{d-1} < y < -1+\frac{i}{d-1}\right\}$$
$$\cup \bigcup_{i=1}^{d-1}\left\{x_i > 0, \quad \frac{i-1}{d-1} < y < \frac{i}{d-1}\right\}.$$

We then perform the same second quantization, as defined per formula (4.2), with $z = \frac{\epsilon(h)}{h^{\frac{1}{2}}}(x_1, \ldots, x_{d-1})$ and $\zeta = \epsilon(h)h^{\frac{1}{2}}(\xi_1, \ldots, \xi_{d-1})$. By the 2-microlocal calculus on $\mathbb{R}^d \times \mathbb{R}^d \times \mathbb{R}^{d-1} \times \mathbb{R}^{d-1}$, the sequence $(u_n)$ still admits a 2-microlocal measure $\mu$. The first two points of Proposition 4.4 are then proved exactly as previously. The proof of the third point is identical since the non-concentration estimates of Proposition 3.1 hold for $d-$dimensional tori. The proof of the microlocal vanishings follow the same lines, and equations (4.15) and (4.16) still hold so that the fifth point is unchanged.

We then conclude in the same manner since Lemma 4.5 is satisfied on any sphere of dimension $2d-3$.

4.3. **Reduction of the general case and conclusion.** We now reduce the general case to the preceding model case. Consider a damping $a$ on the $d$-dimensional torus $\mathbb{T}^d = \mathbb{R}^d/\Gamma$ with $\Gamma$ an orthogonal lattice $A_1\mathbb{Z} \times \cdots \times A_d\mathbb{Z}^d$, $A_1, \ldots, A_d > 0$. Assume that $a$ is a finite sum of characteristic functions of disjoint polyhedrons in $\mathbb{T}^d$, and that all geodesics of the torus which do not enter the interior of the damped region satisfy assumption (1.3).

First we deal with non-closed geodesics. Second, we use an argument from [BZ12, Section 3] to microlocalize near a single closed geodesic, the direction of which is isolated in supp($\nu$). We conclude by using orthonormal changes of coordinates to reduce this geodesic to the one of the model case.

**Proposition 4.6.** — *Non-closed geodesics which are damped in every normal direction (ie satisfying assumption* (1.3)*) intersect the interior of a polyhedron where the damping equals 1.*

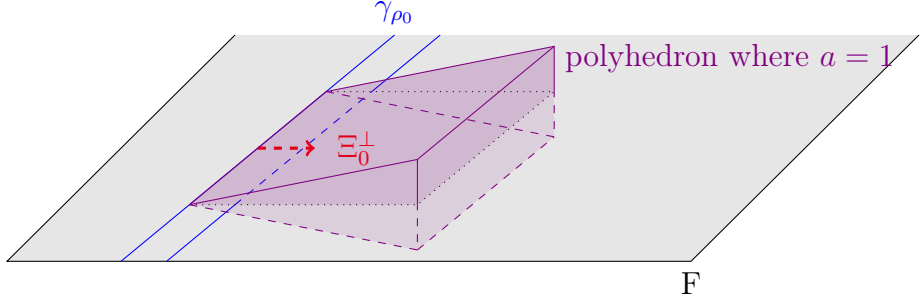

FIGURE 4.2. An illustration of Proposition 4.6. The non-closed geodesic $\gamma_{\rho_0}$ is damped by some polyhedron in the direction $\Xi_0^\perp \in F$. By density in $F$, it enters that polyhedron. The dashed violet edges and paler-shade part of the polyhedron are below the plane $F$. The dashed blue and red lines are inside the polyhedron.

*Proof.* Consider a non-closed geodesic $\gamma_{\rho_0}$ of $\mathbb{T}^d$, $\rho_0 = (X_0, \Xi_0)$. $\gamma_{\rho_0}$ is then dense in some $(X_0 + F)/\Gamma$ where $\rho_0 = (X_0, \Xi_0)$ and $F$ is a subspace of $\mathbb{R}^d$ of dimension strictly greater than one.

The rest of the proof is illustrated in Figure 4.2. Denote $\Xi_0 \in F$ the direction of the geodesic and $\Xi_0^\perp$ some unitary vector in $F$, orthogonal to $\Xi_0$. By (1.3), there exists an interval $I \subset \mathbb{R}$ and some positive $\delta_0$ such that $a = 1$ in an open neighborhood of $\left\{\gamma(s) + \delta \Xi_0^\perp, s \in I, \delta \in \left[\frac{\delta_0}{2}, \delta_0\right]\right\} \subset X_0 + F$. Denote by $V$ such an open neighborhood, then $V \subset \text{Int}(\text{supp}(a))$ and by density of the geodesic, $\gamma_{\rho_0}$ enters $V$. $\square$

The support of $\nu$ thus contains only points $(X, \Xi)$ belonging to closed geodesics, with rational $\Xi \in \mathbb{S}^{d-1}$.

Defining $\pi_\Xi : (X, \Xi) \mapsto \Xi$, we get that $\pi_\Xi(\text{supp}(\nu))$ is countable and closed (using compactness of the torus for the latter). By Baire's lemma, any closed subset with no isolated points (also known as a perfect set) in a complete metric space is uncountable. Thus, $\pi_\Xi(\text{supp}(\nu))$ contains an isolated direction, which we denote $\Xi_0$.

Let us microlocalize around this isolated direction. Consider a closed geodesic $\gamma$ in the support of $\nu$, with direction $\Xi_0$. Consider a neighborhood of $\Xi_0$ which contains no other directions in $\pi_\Xi(\text{supp}(\nu))$, and $\chi(hD_X)$ a Fourier multiplier with a symbol supported in this neighborhood in the $\Xi$ variable. Define

$$v_n = \chi(hD_X)u_n$$

and denote $\nu_v$ the corresponding 1-microlocal measure. $\nu_v$ satisfies

$$(4.19) \qquad \nu_v = \chi(\xi)^2 \nu, \quad \nu_v(T^*\mathbb{T}^d) \geq \nu(\gamma) > 0, \quad (h_n^2 \Delta + 1)v_n = o(h_n),$$

since $h_n^2 \Delta$ and $\chi(hD_x)$ commute. We have thus microlocalized around a single geodesic of $\text{supp}(\nu)$ with an isolated direction.

Since the geodesic is closed, there exist some integers $n_1, \ldots, n_d \in \mathbb{Z}$ such that

$$\Xi_0 = \frac{(n_1 A_1, \ldots, n_d A_d)}{\left(\sum_{j=1}^d n_j^2 A_j^2\right)^{\frac{1}{2}}}.$$

The next lemma explains the induction step which allows to reduce the geodesic to one living in a $(d-1)$-dimensional torus. It is an adaptation of [BZ12, Lemma 2.2], to the $d$-dimensional setting (see also [BG20, Lemma 4.1]).

**Lemma 4.7.** − *Assume that $n_d$ is non-zero. Consider $u_n$ a $\Gamma$-periodic function and denote $\Xi_j$ the $j$-th vector of the canonical basis of $\mathbb{R}^d$ for every $j \in \{1, \ldots, d-2\}$. Denote*

$$\text{(4.20)} \quad \Xi_{d-1} = \frac{(0, \ldots, 0, n_{d-1}A_{d-1}, n_d A_d)}{S_{d-1}}, \Xi_d = \frac{(0, \ldots, 0, -n_d A_d, n_{d-1}A_{d-1})}{S_{d-1}},$$

$$S_{d-1} = \left(n_{d-1}^2 A_{d-1}^2 + n_d^2 A_d^2\right)^{\frac{1}{2}},$$

*so that $(\Xi_j)_{j=1,\ldots,d}$ is a direct orthonormal basis and $\Xi_0 = \sum_{j=1}^{d-2} n_j A_j \Xi_j + S_{d-1} \Xi_{d-1}$. Set*

$$F_{d-1}(x) = \sum_{j=1}^{d} x_j \Xi_j.$$

*Then*

$$\text{(4.21)} \quad u \circ F_{d-1}(x_1 + k_1 A_1, \ldots, x_{d-2} + k_{d-2} A_{d-2}, x_{d-1} + k_{d-1} S_{d-1}, x_d + k_d \alpha) =$$
$$u \circ F_{d-1}(x_1, \ldots, x_{d-2}, x_{d-1} - k_d \beta, x_d), \forall k \in \mathbb{Z}^d, x \in \mathbb{R}^d,$$

*where, for any fixed $p, q \in \mathbb{Z}$,*

$$\alpha = \frac{q n_{d-1} A_{d-1} A_d - p n_d A_d A_{d-1}}{S_{d-1}}, \quad \beta = \frac{p n_{d-1} A_{d-1}^2 + q n_d A_d^2}{S_{d-1}}, \quad \alpha \neq 0.$$

*Proof.* Since $F$ preserves the first $d-2$ vectors of the canonical basis and $u$ is $\Gamma$-periodic, it is sufficient to find some $\alpha, \beta$ such that for any $k_{d-1}, k_d \in \mathbb{Z}$ there exist $p, q \in \mathbb{Z}$ satisfying

$$\text{(4.22)} \quad k_{d-1} S_{d-1} \Xi_{d-1} + k_d \alpha \Xi_d + k_d \beta \Xi_{d-1} = \begin{pmatrix} 0 \\ \vdots \\ 0 \\ pA_{d-1} \\ qA_d \end{pmatrix}$$

The $S_{d-1}$ coefficient in the first term is chosen such that it is sufficient to deal with the case $k_{d-1} = 0$, $k_d = 1$. The last two lines of (4.22) then give an invertible $2 \times 2$ linear system which we solve by taking the scalar product with $\Xi_{d-1}$ and $\Xi_d$. □

We thereby get an orthonormal basis $(\Xi_1, \ldots, \Xi_d)$ and a corresponding coordinate system $(x_1, \ldots, x_d)$ such that $u$ is $(A_1 \mathbb{Z} \times \cdots \times A_{d-2} \mathbb{Z} \times S_{d-1} \mathbb{Z})$-periodic wrt to the $d-1$ first variables, and the geodesic is contained in a $(d-1)$-dimensional torus. Note that if the coefficient $n_d$ is zero, the latter already occurs, so the assumption in the lemma is not restrictive.

Since the non-concentration estimate requires periodicity of $u_n$ in $d-1$ orthogonal directions, we can apply it to $u \circ F_{d-1}$. We get

$$\|u_n \circ F_{d-1}\|_{L^2(\{|x_d| \leq h_n \epsilon(h_n)^{-2}\})} \leq C\epsilon(h_n)^{\frac{1}{2}}.$$

Thus

$$\text{(4.23)} \quad \|u_n \circ F_{d-1}\|_{L^2(F_{d-1}^{-1}(\gamma + B(0, h_n \epsilon(h_n)^{-2})))} \leq C\epsilon(h_n)^{\frac{1}{2}},$$

as a cylinder surrounding the closed geodesic is contained in any $(d-1)$-dimensional slice of identical width.

We now iterate the construction of Lemma 4.7 and show that it preserves the non-concentration estimate. We use the periodicity of $u_n$ with respect to the $(d-1)$ first coordinates to apply the lemma in successive steps: we get $(d-2)$ orthonormal changes of coordinates $F_{d-2}, \ldots, F_1$ (possibly equal to identity if no change of coordinates is necessary) such that for every $j \in \{1, \ldots, d-1\}$:

- each $F_j$ preserves the first $(j-1)$ vectors of the canonical basis,
- $u_n \circ F_{d-1} \circ \ldots \circ F_j$ is periodic wrt its $j$ first variables,
- and the image $F_j^{-1} \circ \ldots \circ F_{d-1}^{-1}(\gamma)$ of the geodesic is contained in a torus of dimension $j$.

Denote $F = F_{d-1} \circ \ldots \circ F_1$. Since every change of variables is isometric, (4.23) becomes

$$(4.24) \qquad \|u_n \circ F\|_{L^2(F^{-1}(\gamma+B(0,h_n\epsilon^{-2}(h_n))))} \leq C\epsilon(h_n)^{\frac{1}{2}}.$$

Since $F^{-1}(\gamma + B(0, h_n\epsilon^{-2}(h_n))) = F^{-1}(\gamma) + B(0, h_n\epsilon^{-2}(h_n))$, the estimate (4.24) is exactly the one that allows to show that the 2-microlocal measure is supported at infinity (see the proof of Proposition 4.4, point 3.). Besides, $F^{-1}(\gamma)$ is a closed geodesic depending on only one cartesian coordinate and satisfying assumption (1.3), so that we are microlocally reduced to the case of subsection 4.2. The rest of the contradiction argument follows the same lines as previously, which completes the proof of Theorem 4.

4.4. **$\zeta.\partial_z$ vector field over the sphere at infinity.** There remains to study the dynamics of vector field $\zeta.\partial_z$ over the sphere at infinity $\mathbb{S}^{2d-3}_\infty$ of $\mathbb{R}^{d-1}_z \times \mathbb{R}^{d-1}_\zeta$, $d \geq 3$, in order to prove Lemma 4.5. The proof relies on symplectic transformations which simplify the expression of the dynamics, and on passing to spherical coordinates (as is done on $\mathbb{S}^1$ in [BG20, Proposition 3.5, 5)]).

First, we reduce the dynamics over $\mathbb{S}^{2d-3}_\infty$ to dynamics over $\mathbb{S}^3_\infty$ where only one angle coordinate is non-constant. We then show that the dynamics over $\mathbb{S}^{2d-3}_\infty$ can be obtained by projecting onto the sphere at infinity the dynamics of $\zeta.\partial_z$ at finite distance from the origin (a straight line traveled at constant velocity $|\zeta(0)|$).

Denote $(z^{(0)}, \zeta^{(0)})$ an initial point belonging to $\mathbb{S}^{2d-3}_\infty$. Since $\mathbb{S}^{2d-3}$ is a compact manifold, the flow of $\zeta.\partial_z$ is defined for all times. We denote $(z(s), \zeta(s))$ the image of $(z^{(0)}, \zeta^{(0)})$ by the flow of $\zeta.\partial_z$ at time $s$. The first result is the following:

**Proposition 4.8.** — (1) *There exists an angle coordinate system $(\theta_1, \ldots, \theta_{2d-3})$ on the sphere at infinity in which the dynamics of $\zeta.\partial_z$ satisfy:*

$$(4.25) \quad \dot{\theta}_1 = -\cos(\theta_2)\sin^2(\theta_1) \quad ; \quad \theta_2 = \theta_2^{(0)} \quad ; \quad \theta_3 = cst = \pm 1 \quad ; \quad \theta_i = 0, \quad \forall i \geq 4.$$

*In other words, in this coordinate system we have $\zeta.\partial_z = -\cos(\theta_2)\sin^2(\theta_1)\partial_{\theta_1}$ over the sphere at infinity.*

(2) $\theta_1(s) \to 0 \mod \pi$ *as $s$ goes to infinity, so that $(z(s), \zeta(s)) \to \left(\pm \frac{\zeta^{(0)}}{|\zeta^{(0)}|}\right)$ if $\zeta^{(0)} \neq 0$. The $(z^{(0)}, \zeta^{(0)})$ such that $\zeta^{(0)} = 0$ are exactly the fixed points of $\zeta.\partial_z$.*

*Proof.* If $\zeta^{(0)} = 0$, $(z^{(0)}, \zeta^{(0)})$ is a fixed point of the dynamic since we are solving a first order initial value problem with zero initial derivative. Thus, we assume $\zeta^{(0)} \neq 0$. We introduce a cartesian coordinate system $(\tilde{z}, \tilde{\zeta})$ such that $\zeta.\partial_z = \tilde{\zeta}.\partial_{\tilde{z}}, \tilde{\zeta}$ has zero coordinates except for the first one and $\tilde{z}$ has zero coordinates except for the first two. Choosing a rotation $R \in SO(d-1)$ such that

$$\tilde{\zeta}^{(0)} := R^T \zeta^{(0)} = \begin{pmatrix} \tilde{\zeta}_1^{(0)} \\ 0 \\ \vdots \\ 0 \end{pmatrix},$$

we set

(4.26) $$\zeta = R\tilde{\zeta} \quad ; \quad z = R^T \tilde{z}.$$

This transformation leaves the vector field unchanged: $\zeta.\partial_z = \tilde{\zeta}.\partial_{\tilde{z}}$. Then along the trajectory of $\zeta.\partial_z$ starting from $(z^{(0)}, \zeta^{(0)})$, $\tilde{\zeta}_i(s)$ vanishes for every $i \geq 2$ and all $s \in \mathbb{R}$. Thus, $\zeta.\partial_z = \tilde{\zeta}.\partial_{\tilde{z}} = \tilde{\zeta}_1 \partial_{\tilde{z}_1}$ along this trajectory.

Proceeding similarly with the $d-1$ last coordinates of $\tilde{z}^{(0)}$, we can guarantee that $\tilde{z}_3^{(0)} = \ldots = \tilde{z}_d^{(0)} = 0$ without changing the coordinates of $\tilde{\zeta}^{(0)}$. Thus, $\tilde{z}_i(s) = 0$ for all $s \in \mathbb{R}$ and every $i \geq 3$.

We now change to spherical coordinates $(r, (\theta_j)_{1 \leq j \leq 2d-3}) \in \mathbb{R}_+^* \times [0, \pi]^{2d-4} \times \mathbb{S}^1$, allowing $r$ to go to $+\infty$. We set $\tilde{z}_1 = r \cos \theta, \tilde{\zeta}_1 = r \sin \theta_1 \cos \theta_2, \tilde{z}_2 = r \sin \theta_1 \sin \theta_2 \cos \theta_3$, and more generally:

(4.27) $$\forall k \in \{1, \ldots, d-1\}, \quad \tilde{z}_k = r \left( \prod_{i=1}^{2k-2} \sin \theta_i \right) \cos \theta_{2k-1} \quad ;$$

$$\tilde{\zeta}_k = r \left( \prod_{i=1}^{2k-1} \sin \theta_i \right) \cos \theta_{2k} \text{ if } k \leq d-2 \quad ; \quad \tilde{\zeta}_{d-1} = r \left( \prod_{i=1}^{2d-3} \sin \theta_i \right).$$

Since all $(\tilde{\zeta}_j)_{2 \leq j \leq d-1}$ and $(\tilde{z}_j)_{3 \leq j \leq d-1}$ vanish, we are reduced to a situation where $\sin \theta_3 = 0$ and $\cos \theta_3 = \pm 1$. We now compute $\tilde{\zeta}_1 \partial_{\tilde{z}_1}$ in the $(r, \theta)$ variables. Since none of the angles $\theta_i$ depend on $\tilde{z}_1$ except for $\theta_1$, we get

$$\tilde{\zeta}_1 \partial_{\tilde{z}_1} = r \sin \theta_1 \cos \theta_2 \left( \frac{\partial r}{\partial \tilde{z}_1} \partial_r + \frac{\partial \theta_1}{\partial \tilde{z}_1} \right).$$

We compute both partial derivatives: $\cos \theta_1 = \frac{\tilde{z}_1}{(|\tilde{z}|^2 + |\tilde{\zeta}|^2)^{\frac{1}{2}}}$ gives $\frac{\partial \theta_1}{\partial \tilde{z}_1} = -\frac{1}{r} \sin \theta_1$. $\frac{\partial r}{\partial \tilde{z}_1} = \frac{\tilde{z}_1}{r}$ is homogeneous of order zero wrt $r$. Given a symbol $q$ that is polyhomogeneous of order zero, $r \partial_r q$ is polyhomogeneous of order -1 so that the $\partial_r$ term above vanishes as $r \to +\infty$. Thus, over the sphere at infinity, we have

$$\forall q \in S^0, \quad \zeta.\partial_z q|_{\mathbb{S}^{2d-3}} = -\cos \theta_2 \sin^2 \theta_1 \partial_{\theta_1} \tilde{q}(X, \Xi, r \to +\infty, \theta)$$

where $\tilde{q}$ is a notation for $q$ in $(X, \Xi, r, \theta)$ coordinates. $\square$

We now give an interpretation of the dynamics of this vector field as the projection onto the sphere at infinity of the flow map of $\zeta.\partial_z$ starting from a point $(z^{(0)}, \zeta^{(0)}) \in \mathbb{S}^1$. This result is illustrated in Figure 4.3.

In the following statement, we denote indifferently points belonging to the unit sphere and the sphere at infinity as they are identified by projection.

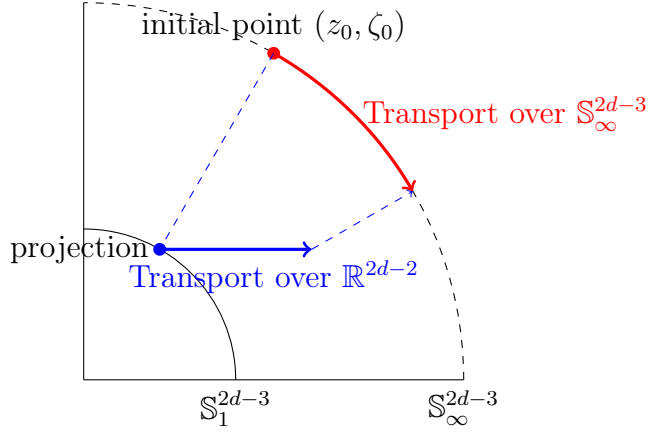

FIGURE 4.3. Illustration of Proposition 4.9 about the transport by vector field $\zeta.\partial_z$ over the sphere at infinity.

**Proposition 4.9.** — *Denote for all $s \in \mathbb{R}$*

$$(4.28) \qquad (\check{z}(s), \check{\zeta}(s)) := \frac{1}{\left(|\tilde{z}^{(0)} + s\tilde{\zeta}^{(0)}|^2 + |\tilde{\zeta}^{(0)}|^2\right)^{\frac{1}{2}}} \left(\tilde{z}^{(0)} + s\tilde{\zeta}^{(0)}, \tilde{\zeta}^{(0)}\right) \in \mathbb{S}^{2d-3},$$

*then*

$$(4.29) \qquad (\tilde{z}(s), \tilde{\zeta}(s)) = (\check{z}(s), \check{\zeta}(s))$$

*where $(\tilde{z}(s), \tilde{\zeta}(s)) \in \mathbb{S}^{2d-3}$ were introduced in (4.26).*

*Proof.* Recall that $(\tilde{z}^{(0)}, \tilde{\zeta}^{(0)})$ has zero-coordinates except maybe for $\tilde{z}_1^{(0)}, \tilde{\zeta}_1^{(0)}, \tilde{z}_2^{(0)}$. Then by (4.28), every coordinate of $(\check{z}, \check{\zeta})$ vanishes except maybe $\check{z}_1, \check{\zeta}_1, \check{z}_2$.

The result then follows from the same passage to angle coordinates as in (4.27). Like before, the vanishing coordinates impose constraints on $\check{\theta}_3$, so that we set

$$\check{z}_1 = r\cos\check{\theta}_1 \quad ; \quad \check{\zeta}_1 = r\sin\check{\theta}_1\cos\check{\theta}_2 \quad ; \quad \check{z}_2 = r\sin\check{\theta}_1\sin\check{\theta}_2\cos\check{\theta}_3 = \pm r\sin\check{\theta}_1\sin\check{\theta}_2.$$

where $\check{\theta}_1, \check{\theta}_2 \in [0, \pi]$ and $\check{\theta}_3 \in [0, \pi]$ if $d > 3$ and $\mathbb{S}^1$ if $d = 3$. Conversely, the $\check{\theta}_j$s satisfy:

$$(4.30) \quad \check{\theta}_1 = \arccos\left(\frac{\check{z}_1}{\sqrt{\sum_{j=1}^p \check{z}_j^2 + \sum_{j=1}^p \check{\zeta}_j^2}}\right), \quad \check{\theta}_2 = \arccos\left(\frac{\check{\zeta}_1}{\sqrt{\sum_{j=2}^p \check{z}_j^2 + \sum_{j=1}^p \check{\zeta}_j^2}}\right).$$

$$(4.31) \quad \begin{aligned} &\text{If } p > 2: \quad \check{\theta}_3 = \arccos\left(\frac{\check{z}_2}{\sqrt{\sum_{j=2}^p \check{z}_j^2 + \sum_{j=2}^p \check{\zeta}_j^2}}\right) = 0 \mod \pi. \\ &\text{If } p = 2: \quad \check{\theta}_3 = \arccos\left(\frac{\check{z}_2}{\sqrt{\sum_{j=2}^p \check{z}_j^2 + \sum_{j=2}^p \check{\zeta}_j^2}}\right) \text{ if } \check{\zeta}_2 \geq 0, \\ &\qquad\qquad\qquad\quad 2\pi - \arccos\left(\frac{\check{z}_2}{\sqrt{\sum_{j=2}^p \check{z}_j^2 + \sum_{j=2}^p \check{\zeta}_j^2}}\right) \text{ otherwise.} \end{aligned}$$

Only the first coordinate of $\tilde{\zeta}^{(0)}$ is non-zero so by (4.28), $\check{\theta}_2$ and $\check{\theta}_3$ are time-independent. Since the $\theta_j$s and $\check{\theta}_j$s have identical initial values, we need to show that

$$\frac{d}{dt}\check{\theta}_1(t) = -\cos\check{\theta}_2 \sin^2\check{\theta}_1(t)$$

to conclude. We differentiate the expression of $\cos(\check{\theta}_1(t))$ to get

$$-\sin(\check{\theta}_1(t))\frac{d\check{\theta}_1(t)}{dt}(t) = \frac{(|\tilde{z}^{(0)} + t\tilde{\zeta}^{(0)}|^2 + |\tilde{\zeta}^{(0)}|^2)\tilde{\zeta}_1^{(0)} - (\tilde{z}_1^{(0)} + t\tilde{\zeta}_1^{(0)})(\tilde{\zeta}^{(0)}.(\tilde{z}^{(0)} + t\tilde{\zeta}^{(0)}))}{(|\tilde{z}^{(0)} + t\tilde{\zeta}^{(0)}|^2 + |\tilde{\zeta}^{(0)}|^2)^{\frac{3}{2}}}.$$

Simplifying this expression by removing vanishing coordinates and using the expression of $\cos(\check{\theta}_2)$, we get

$$-\sin(\check{\theta}_1(t))\frac{d\check{\theta}_1(t)}{dt}(t) = \cos\check{\theta}_2(t)\left(\frac{(\tilde{\zeta}_1^{(0)})^2 + \sum_{j=2}^p (\tilde{z}_j^{(0)})^2}{|\tilde{z}^{(0)} + t\tilde{\zeta}^{(0)}|^2 + |\tilde{\zeta}^{(0)}|^2}\right)^{\frac{3}{2}}$$

$$= \cos\check{\theta}_2(t)(1 - \cos^2\check{\theta}_1(t))^{\frac{3}{2}} = \cos\check{\theta}_2(t)\sin^3\check{\theta}_1(t),$$

hence the result. □

## 5. Proof of Theorem 5

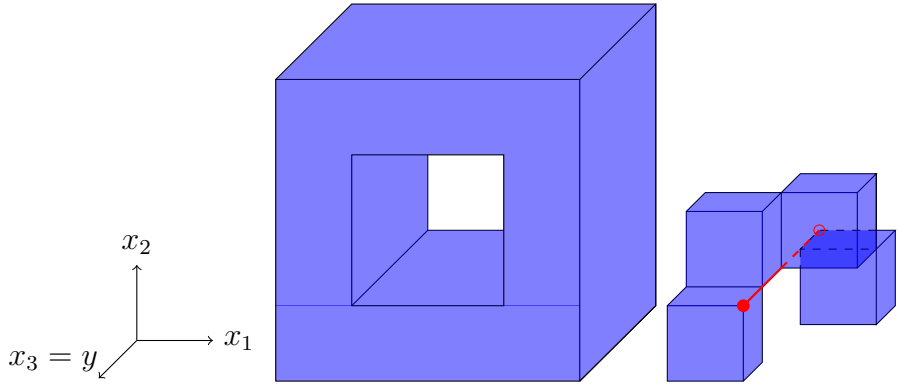

FIGURE 5.1. The damping $a$ is equal to 1 on the blue zone. The flat torus $\mathbb{T}^3$ is identified to the unit cube: $\mathbb{T}^3 \simeq [-1,1]^3/(X \sim X + 2e_i)$. The four cubes on the right fit into the hollow part of the bigger cube. The coordinate system is still $X = (x_1, x_2, y) = (x, y), \Xi = (\xi_1, \xi_2, \eta) = (\xi, \eta)$.

As pointed out in the introduction, the sufficient condition of Theorem 4 is also necessary in dimension 2 [BG20]. The purpose of the present section is to investigate how this condition can be weakened on higher-dimensional tori by proving Theorem 5. We assume that the damping $a$ equals 1 on the blue zone pictured in Figure 5.1. Compared to the model case of section 4.2, the four prisms in the tunnel are reduced to cubes. Said otherwise, we take $\alpha_L = \alpha_R = \alpha_T = \alpha_B = 0$ in (4.8) so that $a$ is the characteristic of the set given in (1.5).

Two main differences emerge in this example compared to the model case of Figure 4.1. First, the bicharacteristics which do not intersect $\text{Int}(\text{supp}(a))$ are the red bicharacteristic $\{x = 0\}$, and parallel ones which raze the faces of the cubes and are contained in the same horizontal and vertical planes as the red one, or more explicitly

$$(5.1) \qquad \{(x_1 = 0, |x_2| \leq 1/2) \text{ or } (|x_1| \leq 1/2, x_2 = 0), \quad \xi = 0, \quad \eta = \pm 1\}.$$

Second, none of these geodesics (in particular not the red one) satisfy Assumption (1.3). The geodesic $\{x = 0\}$ has four undamped directions while the others have two.

In section 5.1 we tackle the first issue. By dropping the $(z_1, \zeta_1)$ or $(z_2, \zeta_2)$ variable in the second quantization of section 4.2, we construct 2-microlocal measures $\mu_1, \mu_2$ in each direction $x_1, x_2$, such that the second microlocal measure $\mu$ is absolutely continuous wrt each. We then show that $\mu$ has no support on any geodesic except the red one. In section 5.2, we perform a third microlocalization around the red geodesic to show that the sequence $u_n$ of quasimodes cannot concentrate around it either. We conclude by discussing the scope of this example and some obstacles to generalizing it.

## 5.1. Directional second microlocalizations.

For a given $j \in \{1, 2\}$, the second quantization in the direction $x_j$ of a symbol $b_j(X, \Xi, z_j, \zeta_j)$ is defined by

$$(5.2) \quad \operatorname{Op}_h^{(\epsilon,j)}(b_j)u(X) = \frac{1}{(2\pi)^3} \int_{\mathbb{R}_Y^3 \times \mathbb{R}_\Xi^3} e^{i(X-Y).\Xi} b_j\left(X, h\Xi, \frac{\epsilon(h)}{h^{\frac{1}{2}}} x_j, \epsilon(h) h^{\frac{1}{2}} \xi_j\right) u(Y) dY d\Xi.$$

This quantization is a specific case of (4.2) with $d = 3$, $q = 1$, so $u_n$ admits two microlocal measures $\mu_1$, $\mu_2$ associated to each quantization $\operatorname{Op}_h^{(\epsilon,1)}$, $\operatorname{Op}_h^{(\epsilon,2)}$. $\mu_1$ (resp. $\mu_2$) measures concentration in the horizontal (resp. vertical) direction, at the 2-microlocal scale $h^{\frac{1}{2}}\epsilon(h)^{-1}$. We also define the second microlocal measure $\mu$ as formerly. The properties of the $\mu_j$s and of $\mu$ are summarized in the following statement:

**Proposition 5.1** (Properties of $\mu_1$ and $\mu_2$). – *The directional 2–microlocal measures have the following properties:*

(1) *Assume that $(h_n^2 \Delta + 1)u_n = \mathcal{O}_{L^2}(1)$, then $\mu_1$, $\mu_2$ and $\mu$ are probability measures.*
(2) *Assume that $(h_n^2 \Delta + 1)u_n = \mathcal{O}_{L^2}(h_n)$ and $\|au_n\|_{L^2} = o(1)$, and define the projection $\pi_j(X, \Xi, z, \zeta) = (X, \Xi, z_j, \zeta_j)$ for $j = 1, 2$. Then, $\mu_j$ is the pushforward measure of $\mu$ by the projection $\pi_j$. Besides, the pushforward measure of each measure $\mu$, $\mu_1$, $\mu_2$ by the projection onto the $(X, \Xi)$ variables is the first microlocal measure $\nu$.*

Thus,

$$(5.3) \quad \operatorname{supp}(\mu) \subset$$
$$\{(X, \Xi, z, \zeta), (x_1 = 0, |x_2| \leq 1/2) \text{ or } (|x_1| \leq 1/2, x_2 = 0), \xi = 0, \eta = \pm 1\}$$

*and the same inclusions hold for $\operatorname{supp}(\mu_j), j = 1, 2$, replacing $(z, \zeta)$ with $(z_j, \zeta_j)$.*
(3) *Assume now that $(h_n^2 \Delta + 1)u_n = \mathcal{O}(h_n \epsilon(h_n))_{L^2}$, then each measure $\mu_j$ is supported in the sphere at infinity $\mathbb{S}_\infty^1$ in the $(z_j, \zeta_j)$ variables.*
(4) *Along the geodesics that raze one of the faces of the four damped cubes in Figure 5.1, directional 2-microlocal measures vanish 2-microlocally in the direction of* *the cube. More precisely, $\mu_1$ satisfies:*

(5.4)
$$\mu_1 \left\{ x_1 = 0, x_2 \in \left[-\frac{1}{2}, 0\right), y \in \left(-1, -\frac{1}{2}\right), \xi = 0, \eta = \pm 1, z_1 > 0 \right\} = 0,$$
$$\mu_1 \left\{ x_1 = 0, x_2 \in \left(0, \frac{1}{2}\right], y \in \left(-\frac{1}{2}, 0\right), \xi = 0, \eta = \pm 1, z_1 > 0 \right\} = 0,$$
$$\mu_1 \left\{ x_1 = 0, x_2 \in \left(0, \frac{1}{2}\right], y \in \left(0, \frac{1}{2}\right), \xi = 0, \eta = \pm 1, z_1 < 0 \right\} = 0,$$
$$\mu_1 \left\{ x_1 = 0, x_2 \in \left[-\frac{1}{2}, 0\right), y \in \left(\frac{1}{2}, 1\right), \xi = 0, \eta = \pm 1, z_1 < 0 \right\} = 0,$$

*and similar vanishings hold for $\mu_2$, which we do not state for the sake of conciseness.*

(5) *Consider $j \in \{1, 2\}$, then the microlocal measure $\mu_j$ satisfies the conservation law*

(5.5) $$(\eta \partial_y + \zeta_j \partial_{z_j}) \mu_j = 0.$$

*Since $\mu_j$ is supported in the sphere at infinity in the $(z_j, \zeta_j)$ variables, we can use polar coordinates $(z_j, \zeta_j) = (r \cos \theta, r \sin \theta)$, $r \to +\infty$ with $\theta \in \mathbb{R}/2\pi\mathbb{Z}$. In these coordinates, equation (5.5) becomes*

(5.6) $$(\eta \partial_y - \sin^2(\theta) \partial_\theta) \mu_j = 0.$$

(6) *The two preceding items imply that*

(5.7)
$$\operatorname{supp}(\mu_1) \subset \{(X, \Xi, z, \zeta) : |x_1| \leq \frac{1}{2}, x_2 = 0, \xi = 0, \eta = \pm 1, (z_1, \zeta_1) \in \mathbb{S}^1_\infty\}$$
$$\operatorname{supp}(\mu_2) \subset \{(X, \Xi, z, \zeta) : x_1 = 0, |x_2| \leq \frac{1}{2}, \xi = 0, \eta = \pm 1, (z_2, \zeta_2) \in \mathbb{S}^1_\infty\}$$

*Since the support of $\mu$ is a subset of those of $\mu_1$ and $\mu_2$, we get that $\mu$ is supported in the undamped great circles of the sphere at infinity $\mathbb{S}^3_\infty$ along the red bicharacteristics:*

(5.8) $$\operatorname{supp}(\mu) \subset \{x = 0, \xi = 0, \eta = \pm 1, (z_1, \zeta_1) = (0, 0) \text{ or } (z_2, \zeta_2) = (0, 0)\}.$$

*Proof.* The first point derives from $h_n$-oscillation.

As in the proof of the second point of Proposition 4.4, we have that

$$\langle \mu_j, b \rangle = \langle \mu, b \circ \pi_j \rangle$$

for any symbol $b(X, \Xi, z_j, \zeta_j)$. The rest of the proof of the second point is identical.

Concerning (3), we consider a function $\chi \in C_c^\infty(T^*\mathbb{T}^3 \times \mathbb{R}^2_{z_j, \zeta_j})$. There exists some $A > 0$ such that $\chi$ is supported in the set $\{|x_j| \leq A\}$. Then by Proposition 3.1,

(5.9)
$$|\langle \operatorname{Op}^{(\epsilon, j)}_{h_n}(\chi) u_n, u_n \rangle_{L^2}| = |\langle \operatorname{Op}^{(\epsilon, j)}_{h_n}(\chi) u_n, \mathbb{1}_{\{|x_j| \leq A h_n^{\frac{1}{2}} \epsilon^{-1}(h_n)\}} u_n \rangle_{L^2}|$$
$$\leq \|\operatorname{Op}^{(\epsilon, j)}_{h_n}(\chi) u_n\|_{L^2} \|u_n\|_{L^2(\{|x_j| \leq A h_n^{\frac{1}{2}} \epsilon^{-1}(h_n)\})}$$
$$\leq C \epsilon^{\frac{1}{2}}(h_n) \to_{n \to +\infty} 0.$$

We now prove the first microlocal vanishing of $\mu_1$. All other vanishings are showed similarly. Hence we are set to show that

$$\mu_1 \left\{ x_1 = 0, x_2 \in \left[-\frac{1}{2}, 0\right), y \in \left(-1, -\frac{1}{2}\right), z_1 > 0 \right\} = 0.$$

Consider some positive $\delta$ and $\delta_0$ and the symbol

$$(5.10) \quad b_1(X, \Xi, z_1, \zeta_1) = \chi_1\left(\frac{x_1}{\delta}\right) \chi_2(x_2)\chi_3(y)\chi_1(\xi)\chi_1(\eta - 1)\psi\left(\frac{z_1}{\delta_0|\zeta_1|}\right)\psi(|z_1|^2 + |\zeta_1|^2)$$

where $\chi_i$, $i = 1, 2, 3$ and $\psi$ are $[0, 1]$-valued cutoff functions in $C^\infty(\mathbb{R})$ such that:
- $\chi_1$ equals 1 over $[-\frac{1}{2}, \frac{1}{2}]$ and is supported in $[-1, 1]$,
- $\chi_2$ equals 1 over $[-\frac{1}{2}, \delta_0]$ and is supported in $[-\frac{1}{2} - \delta_0, 0]$,
- $\chi_3$ equals 1 over $[-1 + \delta_0, -\frac{1}{2} - \delta_0]$ and is supported in $[-1, -\frac{1}{2}]$.
- $\psi$ equals 1 over $[1, +\infty)$ and is supported inside $[\frac{1}{2}, +\infty)$.

We then proceed as in Proposition 4.4 to prove that $\langle \mu_1, b_1 \rangle = 0$ and deduce the 2-microlocal vanishing from the limit $\delta_0 \to 0$, giving the fourth point.

In the fifth point, (5.5) follows from computing the bracket $[h_n^2\Delta + 1, \operatorname{Op}_{h_n}^{(\epsilon,j)}(q)]$ and passing to the limit $h \to 0$. Indeed, (4.15) becomes

$$(5.11) \quad \frac{1}{2ih_n}[h_n^2\Delta + 1, \operatorname{Op}_{h_n}^{(\epsilon,j)}(q)] = \frac{1}{2i}[h_n\Delta, \operatorname{Op}_{h_n}^{(\epsilon,j)}(q)]$$
$$= \operatorname{Op}_{h_n}^{(\epsilon,j)}((\xi.\partial_x + \eta\partial_y + \zeta_j\partial_{z_j})q) - i\frac{h_n}{2}\operatorname{Op}_{h_n}^{(\epsilon,j)}(\Delta_X q) - i\frac{h_n}{2}(\epsilon(h_n)h_n^{-\frac{1}{2}})\operatorname{Op}_{h_n}^{(\epsilon,j)}((\partial_{x_j}.\partial_{z_j})q)$$
$$- i\frac{h_n}{2}(\epsilon(h_n)h_n^{-\frac{1}{2}})^2\operatorname{Op}_{h_n}^{(\epsilon,j)}(\partial_{z_j}^2 q).$$

(5.5) then derives from the same argument as in the proof of (4.11). Thus for any polyhomogeneous symbol $b_j$ of degree 0,

$$(5.12) \quad \begin{aligned} \zeta_j \partial_{z_j} q_j &= \lim_{r \to +\infty} (-\sin^2(\theta)\partial_\theta + r\cos\theta\sin\theta\partial_r)\tilde{q}(X, \Xi, r, \theta) \\ &= -\sin^2(\theta)\partial_\theta \lim_{r \to +\infty} \tilde{q}(X, \Xi, r, \theta). \end{aligned}$$

Note that this computation from [BG20, p. 640] is the equivalent on $\mathbb{S}^1$ of Proposition 4.8 and Lemma 4.5.

Lastly, we prove the first inclusion of (5.7) by contradiction, as was done in subsection 4.2. Assume that $\operatorname{supp}(\mu_1)$ contains a point $(x, y, \xi = 0, \eta = \pm 1, \theta)$ with $x_1 = 0$, $x_2 \in [-\frac{1}{2}, \frac{1}{2}] \setminus \{0\}$ and $\theta$ defined by $(z_1, \zeta_1) = r(\cos\theta, \sin\theta)$, $r \to +\infty$. We denote $\phi_s(\theta_0)$ the flow of $\dot\theta = -\sin^2\theta$ at time $s$ starting from a point $\theta_0$. The conservation law for $\mu_1$ gives that

$$(5.13) \quad (x, y + s\eta, \xi = 0, \eta, \phi_s(\theta))$$

also belongs to $\operatorname{supp}(\mu_1)$ for any $s \in \mathbb{R}$. A quick study of the flow of equation $\dot\theta = -\sin^2\theta$ (see [BG20, p. 641]) shows that $\phi_s(\theta_0)$ converges to $0 \pmod \pi$ so that the corresponding $(z_1(s), \zeta_1(s))$ converges to $(\pm 1, 0)$ on the sphere at infinity as $s$ goes to $+\infty$. For $s$ large enough, this means that a point of the form (5.13) belongs to one of the zones where $\mu$ vanishes two-microlocally, which gives a contradiction.

The proof for $\operatorname{supp}(\mu_2)$ is identical. We come to the localization of $\operatorname{supp}(\mu)$. Like in Section 4.2, $\mu$ satisfies the conservation law (4.11) and the 2-microlocal vanishings (4.10) with $\alpha_R = \alpha_T = \alpha_L = \alpha_B = 0$. By Lemma 4.5 and the contradiction argument of section 4.2, $\mu$ is supported in great circles of the sphere at infinity $\mathbb{S}^3_{z,\zeta}$ that do not intersect any of the four quadrants $\{\varepsilon_1 z_1 > 0, \varepsilon_2 z_2 > 0\}$, $\varepsilon_1, \varepsilon_2 \in \{-1, 1\}$. This implies inclusion (5.8). □

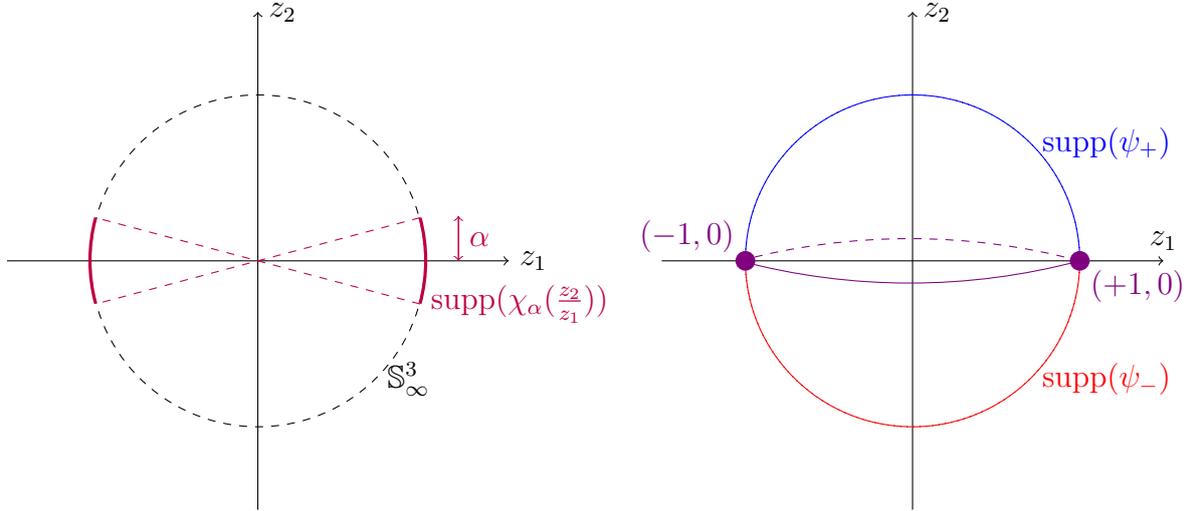

FIGURE 5.2. On the left-hand side, we show the $|z_2| \ll |z_1|$ zone of the sphere at infinity where the $\chi_\alpha(\frac{z_2}{z_1})$ cutoff localizes. On the right-hand-side, we display the supports of the $\psi_+$ and $\psi_-$ cutoffs over the sphere at infinity (in the $(z_1, z_2)$ variables). The violet circle in the horizontal plane represents the $\{z_2 = 0, \zeta_2 = 0\}$ great circle of $\mathbb{S}^3$ (the third coordinate is then $\zeta_1$). The goal of the third microlocalization is to deal with points of $\operatorname{supp}(\mu)$ belonging to this great circle. The flow of $\zeta.\partial_z$ transports any point in this great circle close to one of the two violet disks $(z_1, z_2) = (\pm 1, 0)$ of $\mathbb{S}^3_\infty$ as time goes to infinity. The goal of the third microlocalization is to distinguish points $(z_1, z_2) = (\pm 1, 0^+)$ and $(\pm 1, 0^-)$.

5.2. **Third microlocalization near isolated undamped normal directions.** To conclude the contradiction argument and show that the second microlocal measure $\mu$ vanishes everywhere, we need to show that it has no weight over the undamped directions $(z_1, \zeta_1) = (0, 0)$ and $(z_2, \zeta_2) = (0, 0)$ of the sphere at infinity. To do this, we perform a third microlocalization which allows to add inhomogeneity in the symbols considered and split $\mu$ into two 2-microlocal measures $\mu_+, \mu_-$ supported in each hemisphere of the sphere at infinity.

As an example, we focus on the undamped direction $(z_1, z_2, \zeta_1, \zeta_2) = (z_1, 0, \zeta_1, 0)$ with $z_1 > 0$. Under the flow of $\zeta.\partial_z$, $\zeta_1$ becomes arbitrarily small for long times so that we want to localize on the region of the sphere where $z_1 \gg |\zeta|$ and $z_1 \gg |z_2|$. We denote $b_-$ a symbol of $S^0_H$ localizing in this region:

$$(5.14) \quad b_-(X, \Xi, z, \zeta) = \chi_\delta(|x|)\chi_\delta(|\xi|)\tilde{\chi}(y)\chi_\delta(\eta - 1)\psi\left(\frac{z_1}{\alpha|\zeta|}\right)\chi_\alpha\left(\frac{z_2}{z_1}\right)\psi(|z|^2 + |\zeta|^2),$$

where all cut-offs are smooth functions and satisfy:

- $\chi_\delta$ equals 1 over $[-\frac{\delta}{2}, \frac{\delta}{2}]$ and 0 out of $[-\delta, \delta]$. $\chi_\alpha$ is the same with a different small parameter.
- $\tilde{\chi}$ equals 1 over $[-1+\delta_0, -\frac{1}{2} - \delta_0]$ and 0 out of $[-1, -\frac{1}{2}]$.
- $\psi$ equals 1 over $[1, +\infty)$ and 0 over $(-\infty, \frac{1}{2}]$.

$b_-$ belongs to the smooth polyhomegeneous symbol class $S_H^0$ of the 2-microlocal pseudodifferential calculus. The third microlocalization consists in multiplying $\mathrm{Op}_{h_n}^{(\epsilon)}(b_-)$ by $\psi(-\epsilon^{\frac{1}{2}}(h)z_2)$.

The operator obtained in this way localizes in the $\{y \in (-1, -\frac{1}{2})\}$ portion of the red geodesic in the 1-microlocal variables. In the 2-microlocal variables, it localizes near point $(z_1, z_2, \zeta_1, \zeta_2) = (+1, 0, 0, 0)$ and in the $\{z_1 > 0, z_2 < 0\}$ quadrant of the sphere at infinity, hence in a region of the torus where the damping equals 1 (see Figure 5.1). This cutoff is very similar to the cutoff $b_1$ introduced in (5.10), except that here the red geodesic razes an edge of the damped region, whereas the geodesics we dealt with in Section 5.1 razed a face: we therefore need to localize in a quadrant of the sphere at infinity rather than a hemisphere.

Note also that since this third microlocalization procedure deals only with points belonging to the great circle $\{z_2 = 0, \zeta_2 = 0\}$ of $\mathrm{supp}(\mu)$, one can change $\psi\left(\frac{z_1}{\alpha|\zeta|}\right)$ into $\psi\left(\frac{z_1}{\alpha|\zeta_1|}\right)$ in (5.14). Points belonging to $\{z_1 = 0, \zeta_1 = 0\}$ (the other great circle in which $\mathrm{supp}(\mu)$ may be localized) must be dealt with using a similar third microlocalization in the $z_1$ direction.

The procedure is explained in the following statement and illustrated in Figure 5.2:

**Proposition 5.2.** – (1) Denote $\psi_\pm(x_2) = \psi\left(\pm \frac{\epsilon^{\frac{3}{2}}(h)}{h^{\frac{1}{2}}} x_2\right)$. Denote also $\psi_\pm$ the associated multiplication operators. Then there exist two positive Radon measures $\mu_\pm$ such that for every $b \in S_H^0$,

$$\langle \psi_\pm \mathrm{Op}_{h_n}^{(\epsilon)}(b) u_n, u_n \rangle_{L^2} \to_{n \to +\infty} \langle \mu_\pm, b \rangle.$$

(2) For any 2-microlocal symbol $b \in S_H^0$, we have

$$\langle (\psi_+ + \psi_-) \mathrm{Op}_{h_n}^{(\epsilon)}(b) u_n, u_n \rangle_{L^2} \to_{n \to +\infty} \langle \mu, b \rangle,$$

so that $\mu = \mu_+ + \mu_-$. In particular, $\mathrm{supp}(\mu_\pm) \subset \mathrm{supp}(\mu)$.

(3) $\mu_+$ and $\mu_-$ are supported in the sphere at infinity in the $(z, \zeta)$ variables. Besides, $\mu_+$ is supported in the $\{z_2 \geq 0\}$ closed hemisphere of the sphere at infinity and $\mu_-$ is supported in the $\{z_2 \leq 0\}$ hemisphere.

(4) The $\mu_\pm$ measures satisfy the conservation laws

(5.15) $$\mathbb{1}_{\pm z_2 \geq 0} (\eta \partial_y + \zeta . \partial_z) \mu_\pm = 0.$$

Since $\mathrm{supp}(\mu_\pm) \subset \{\pm z_2 \geq 0\}$, the characteristic functions can be omitted so the $\mu_\pm$ satisfy the same conservation law (4.11) as $\mu$.

(5) Consider $b_-$ the symbol given by (5.14), then $\langle \mu_-, b_- \rangle = 0$. Similarly, replacing $\tilde{\chi}$ by a function equal to 1 over $[-\frac{1}{2} + \delta_0, -\delta_0]$ and supported in $[-\frac{1}{2}, 0]$ in (5.14) defines a symbol $b_+$ such that $\langle \mu_+, b_+ \rangle = 0$. Thus for some small enough positive

$$\alpha,$$

(5.16)
$$\mu_- \left\{ x = 0, y \in \left(-1, -\frac{1}{2}\right), \xi = 0, \eta = \pm 1, z_1 > \frac{2}{\alpha}|z_2|, z_1 > \alpha|\zeta| \right\} = 0,$$
$$\mu_+ \left\{ x = 0, y \in \left(-\frac{1}{2}, 0\right), \xi = 0, \eta = \pm 1, z_1 > \frac{2}{\alpha}|z_2|, z_1 > \alpha|\zeta| \right\} = 0.$$

A similar procedure for negative $z_1$ gives

(5.17)
$$\mu_+ \left\{ x = 0, y \in \left(0, \frac{1}{2}\right), \xi = 0, \eta = \pm 1, z_1 < -\frac{2}{\alpha}|z_2|, z_1 < -\alpha|\zeta| \right\} = 0,$$
$$\mu_- \left\{ x = 0, y \in \left(\frac{1}{2}, 1\right), \xi = 0, \eta = \pm 1, z_1 < -\frac{2}{\alpha}|z_2|, z_1 < -\alpha|\zeta| \right\} = 0.$$

We add a few comments concerning this result and the necessity of introducing the cutoff functions $\psi_\pm$. As illustrated in Figure 5.2, these cutoffs localize in the region $\pm z_2 > 0$, respectively. This localization property is key to proving point (5) of Proposition 5.2. Indeed, by Proposition 5.1, the measure $\mu$ must concentrate 2-microlocally near points $(z_1, z_2, \zeta_1, \zeta_2) = (\pm 1, 0, 0, 0)$ or $(0, \pm 1, 0, 0)$ of the sphere at infinity (see Figure 5.2). The vanishing properties (5.16) then show that concentration in these directions is impossible because they are tangent to faces of the damped region. Tangency to faces in the $z_2 > 0$ region is captured using cutoff $\psi_+$ and gives vanishing properties for $\mu_+$, while tangency to faces in the $z_2 < 0$ region gives vanishing properties for $\mu_-$.

The $\frac{\epsilon^{\frac{3}{2}}(h)}{h^{\frac{1}{2}}} x_2$ scaling in $\psi_\pm$ is chosen so that $1 - (\psi_+ + \psi_-)$ is localized in the non-concentration slice of Proposition 3.1, and all the terms in (5.19) vanish except for one. The first fact is a key to proving point (2), and the second one to proving point (4).

*Proof.* The $\psi_\pm$ are multiplication operators by smooth, uniformly bounded functions so they are $h$-uniformly bounded operators over $L^2_{ul}$. $\operatorname{Op}^{(\epsilon)}_{h_n}(b)$ is also a uniformly bounded operator.

Besides, $\psi_\pm \operatorname{Op}^{(\epsilon)}_{h_n}(b) = \operatorname{Op}^{(\epsilon)}_{h_n}(\psi_\pm(\epsilon^{\frac{1}{2}} z_2) b)$ is the 2-microlocal quantization of a function in the generalized symbol class $S^0$ of our 2-microlocal calculus. Since

$$b \geq 0 \Rightarrow \psi(\epsilon(h)^{\frac{1}{2}} z_2) b \geq 0,$$

the Gårding inequality (4.6) still holds. The classical proof of existence of semiclassical measures then gives existence and positivity of the Radon measures $\mu_\pm$, hence the first point.

The second point derives from the non-concentration estimate of Proposition 3.1. The function $1 - \left( \psi\left(\frac{\epsilon^{\frac{3}{2}}(h) x_2}{h^{\frac{1}{2}}}\right) + \psi\left(-\frac{\epsilon^{\frac{3}{2}}(h) x_2}{h^{\frac{1}{2}}}\right) \right)$ is indeed supported in the $\{|x_2| \leq h^{\frac{1}{2}} \epsilon(h)^{-\frac{3}{2}}\}$ slice. As in the previous proof, this yields

$$\langle (1 - (\psi_+ + \psi_-)) \operatorname{Op}^{(\epsilon)}_{h_n}(b) u_n, u_n \rangle \to_{n \to +\infty} 0,$$

hence the result.

Thus, $\mu_\pm \ll \mu$ and the 3-microlocal measures are supported in the sphere at infinity in the $(z, \zeta)$ variables.

To prove the second part of statement (3), consider a symbol $b$ in $S^0_H$ supported in the $\{z_2 > 0\}$ (resp. $\{z_2 < 0\}$) part of the sphere at infinity. Then $\psi_- b$ (resp. $\psi_+ b$) vanishes in a neighborhood of the sphere at infinity, so that $\langle \mu_-, b \rangle = 0$ (resp. $\langle \mu_+, b \rangle = 0$).

We derive point (4) from computing the bracket $\frac{1}{2ih}[h^2\Delta + 1, \psi_\pm \operatorname{Op}_h^{(\epsilon)}(b)]$:

$$\frac{1}{2ih}[h^2\Delta + 1, \psi\left(\pm\frac{\epsilon^{\frac{3}{2}}(h)}{h^{\frac{1}{2}}}x_2\right)\operatorname{Op}_h^{(\epsilon)}(b)]u$$

$$= \frac{h}{2i}\Delta(\psi_\pm \operatorname{Op}_h^{(\epsilon)}(b)u) - \psi_\pm \operatorname{Op}_h^{(\epsilon)}(b)\left(\frac{h}{2i}\Delta u\right)$$

(5.18)
$$= \frac{h}{2i}\left(\frac{\epsilon^3(h)}{h}\psi''_\pm \operatorname{Op}_h^{(\epsilon)}(b)u \pm 2\frac{\epsilon^{\frac{3}{2}}(h)}{h^{\frac{1}{2}}}\psi'_\pm \partial_{x_2}(\operatorname{Op}_h^{(\epsilon)}(b)u) + \psi_\pm[\Delta, \operatorname{Op}_h^{(\epsilon)}(b)]u\right)$$

$$= \frac{1}{2i}\left(\psi''_\pm \epsilon^2(h)\operatorname{Op}_h^{(\epsilon)}(b)u \pm 2\psi'_\pm \epsilon^{\frac{3}{2}}(h)h^{\frac{1}{2}}\operatorname{Op}_h^{(\epsilon)}\left(i\xi_2 b + \partial_{x_2}b + \frac{\epsilon(h)}{h^{\frac{1}{2}}}\partial_{z_2}b\right)u\right.$$

$$\left. + \psi_\pm[h\Delta, \operatorname{Op}_h^{(\epsilon)}(b)]u\right)$$

Computing the last bracket using (4.15) gives

(5.19) $\frac{1}{2ih}[h^2\Delta + 1, \psi\left(\pm\frac{\epsilon^{\frac{3}{2}}(h)}{h^{\frac{1}{2}}}x_2\right)\operatorname{Op}_h^{(\epsilon)}(b)]$

$$= \frac{1}{2i}\psi''_\pm \epsilon^3(h)\operatorname{Op}_h^{(\epsilon)}(b) \pm \frac{1}{i}[\epsilon(h)^{\frac{1}{2}}\psi'_\pm \operatorname{Op}_h^{(\epsilon)}(i\zeta_2 b) + \epsilon^{\frac{3}{2}}(h)h^{\frac{1}{2}}\operatorname{Op}_h^{(\epsilon)}(\partial_{x_2}b) + \epsilon^{\frac{5}{2}}(h)\operatorname{Op}_h^{(\epsilon)}(\partial_{z_2}b)]$$

$$+ \frac{1}{2i}\psi_\pm\left[\operatorname{Op}_h^{(\epsilon)}(2i(\eta\partial_y + \xi.\partial_x + \zeta.\partial_z)b) + h\operatorname{Op}_h^{(\epsilon)}(\Delta_X b)\right.$$

$$\left. + \epsilon(h)h^{\frac{1}{2}}\operatorname{Op}_h^{(\epsilon)}(\partial_x.\partial_z b) + \epsilon(h)^2\operatorname{Op}_h^{(\epsilon)}(\partial_z^2 b)\right].$$

Except for the $\frac{1}{2i}\psi_\pm \operatorname{Op}_h^{(\epsilon)}(2i(\eta\partial_y + \xi.\partial_x + \zeta.\partial_z)b)$ term, all terms vanish as $h$ goes to zero. Over the sphere at infinity, $\psi\left(\pm\frac{\epsilon^{\frac{3}{2}}(h)}{h^{\frac{1}{2}}}x_2\right)$ is equal to $\mathbb{1}_{\pm z_2 \geq 0}$, giving the result.

To prove the last point, fix some small $\alpha > 0$ independent of $h$ and consider the symbol $b_-$ given by (5.14), then for any $n$, $\operatorname{supp}(\psi_- \operatorname{Op}_{h_n}^{(\epsilon)}(b_-)u_n) \subset \{x_1 > 0, x_2 < 0, y \in (-1, -\frac{1}{2})\}$ so that

$$\langle \psi_- \operatorname{Op}_{h_n}^{(\epsilon)}(b_-)u_n, u_n\rangle = \langle \mathbb{1}_{\{x_1>0, x_2<0, y\in(-1,-\frac{1}{2})\}}\psi_- \operatorname{Op}_{h_n}^{(\epsilon)}(b_-)u_n, u_n\rangle$$

$$= \langle \psi_- \operatorname{Op}_{h_n}^{(\epsilon)}(b_-)u_n, \mathbb{1}_{\{x_1>0, x_2<0, y\in(-1,-\frac{1}{2})\}}au_n\rangle$$

Since $\|au_n\|_{L^2} \to 0$ and $\psi_- \operatorname{Op}_{h_n}^{(\epsilon)}(b_-)u_n$ is bounded uniformly wrt small $h_n$, we get that $\langle \mu_-, b_-\rangle = 0$. We conclude like in the proof of Proposition 5.1, point (4). Since $b_-$ is equal to 1 over $B_x(0, \frac{\delta}{2}) \times B_\xi(0, \frac{\delta}{2}) \times (-1 + \delta_0, -\frac{1}{2} - \delta_0) \times \{z_1 > \alpha|\zeta|, |z_2| < \frac{\alpha}{2}z_1\}$, we get that

$$\mu_-\left\{x = 0, \xi = 0, y \in \left(-1 + \delta_0, -\frac{1}{2} - \delta_0\right), z_1 > \alpha|\zeta|, |z_2| < \frac{\alpha}{2}z_1\right\} = 0.$$

Taking the limit $\delta_0 \to 0$ yields the first vanishing property for $\mu_-$. The other proofs are similar. $\square$

Note that the microlocal vanishings which are proved in point (5) of Proposition 5.2 only concern the measures $\mu_\pm$ rather than $\mu_\infty$. As a consequence, the conservation law $(\eta \partial_y + \zeta.\partial_z)\mu_\pm = 0$ must be applied separately to $\mu_+$ and $\mu_-$ (while in subsection 4.2 it was applied to $\mu$ as a whole).

Consider a point $(x = 0, y, \xi = 0, \eta = \pm 1, z_0, \zeta_0)$ which we assume to belong to supp$(\mu)$. By Proposition 5.1, $(z_0, \zeta_0)$ belongs to one of the great circles $\{(z_1, \zeta_1) = (0, 0)\}$ or $\{(z_2, \zeta_2) = (0, 0)\}$ of the sphere at infinity. Assume the latter. Since $\mu = \mu_+ + \mu_-$, this point must belong either to supp$(\mu_+)$ or to supp$(\mu_-)$. By Lemma 4.5, $\phi_s(z_0, \zeta_0)$ converges to some point in $(z_1, z_2, \zeta_1, \zeta_2) = (\pm 1, 0, 0, 0)$ as $s \to +\infty$. However, the microlocal vanishings of Proposition 5.2, (5) show that both $\mu_+$ or $\mu_-$ in a neighborhood of these points somewhere along the bicharacteristic $\{x = 0, \xi = 0, \eta = \pm 1\}$. Like in subsection 4.2, this means that the measure $\mu$ cannot have support in the great circle $\{(z_2, \zeta_2) = (0, 0)\}$.

If $(z_0, \zeta_0)$ belongs to the great circle $\{(z_1, \zeta_1) = (0, 0)\}$, its limit under the transport of the flow $\zeta.\partial_z$ is one of the points $(z_1, z_2, \zeta_1, \zeta_2) = (0, \pm 1, 0, 0)$. In that case, we perform the same third microlocalization in the $|z_2| \gg |z_1|$ region of the sphere at infinity, and replace cutoffs $\psi_\pm$ with multipliers in the $x_1$ variable. The same arguments give that $\mu$ cannot have support in the great circle $\{(z_1, \zeta_1) = (0, 0)\}$, thus $\mu$ is identically zero. This contradicts that $\mu$ is a probability measure, thus proving stabilization for the damping of Figure 5.1.

We finish with some concluding remarks concerning possible generalizations of this counter-example (and some obstacles to them):

- Based on the result and techniques of the present section, we conjecture that a necessary and sufficient condition for uniform stabilization of the wave equation on 3-dimensional tori using damped polyhedrons is for every geodesic to be damped in every normal direction but a finite number. It is yet unclear what a similar condition should be in dimensions higher than 3. Proof of the necessary condition would likely resemble that of Nicolas Burq and Patrick Gérard in [BG20, Section 5], using the necessity of their generalized geometric control condition. The sufficient condition also seems difficult to prove with the techniques of the present article. The main obstacle is that both the directional second microlocalization of section 5.1 and the third microlocalization of section 5.2 can be performed only in planes where the non-concentration estimates hold (this is crucial to Proposition 5.1, point (3) and Proposition 5.2, point (2)).
- For the proof to hold, such planes need to be rational (that is, they cannot be dense in the torus). This is needed due to the non-concentration estimates, which require periodicity of quasimodes in $(d-1)$ directions. Hence, our techniques allow to deal with one-directional sheets of geodesics contained in a rational hyperplane of the torus (as in subsection 5.1) and intersections of such sheets (using the third microlocalization), but not with cases where such sheets of geodesics are contained in irrational planes. Thus, they allow to prove stabilization on $\mathbb{T}^3$ in the case where every geodesic is damped in all normal directions but a finite number of *rational* ones.
- Using as many additional microlocalizations as the dimension, it may even be possible to obtain generalizations of this sufficient condition to higher-dimensional tori as long as the damping is a characteristic function of polyhedrons with rational faces. However, the generalized geometric control condition conjectured by N.

Burq and P. Gérard in [BG20, Section 5] indicates that rationality or irrationality of the undamped normal directions to the geodesics should not be relevant to uniform stabilization. For that reason, we have chosen not to confront the technicalities of such generalizations until additional work allows to bypass the non-concentration estimates.

Laboratoire de Mathématiques d'Orsay, Université Paris-Saclay, Bâtiment 307, 91405 Orsay Cedex

*Email address*: marc.rouveyrol@universite-paris-saclay.fr